\title{Coefficient convexity of divisors of $x^{n}-1$}
\author{Andreas Decker and Pieter Moree}
\documentclass[12pt]{article}
\usepackage{amssymb, latexsym, amsfonts}
\def\@ptsize{2}
\newtheorem{Thm}{Theorem}

\newtheorem{Lem}{Lemma}

\newtheorem{cor}{Corollary}

\newcommand{\qed}{\hfill $\Box$}

\begin{document}
\date{}
\maketitle
{\def\thefootnote{}
}
{\def\thefootnote{}
\begin{abstract}
\noindent We say a polynomial $f\in {\mathbb Z}[x]$ is {\it strongly coefficient convex} if the
set of coefficients of $f$ consists of consecutive integers only. We establish various
results suggesting that the divisors of $x^n-1$ that are in ${\mathbb Z}[x]$ have the
tendency to be strongly coefficient convex and have small coefficients. The case where
$n=p^2q$ with $p$ and $q$ primes is studied in detail.
\end{abstract}
\section{Introduction}
Let $f(x)=\sum_{j=0}^{\infty} c_jx^j$ be a polynomial. We put
${\cal C}_0(f)=\{c_j\}$. Trivially ${\cal C}_0(f)={\cal C}(f)\cup \{0\}$, where
${\cal C}(f)=\{c_j:0\le j\le {\rm deg}(f)\}$ denotes the set of coefficients of $f$.
If there exist integers $a$ and $b$ such that ${\cal C}_0(f)$ consists
of the consecutive integers $a,a+1,\ldots,b-1,b$, then we say that $f$ is {\it coefficient convex}
and write ${\cal C}_0(f)=[a,b]$. If ${\cal C}(f)=[a,b]$, then we say that $f$ is {\it strongly coefficient convex}.
We say that $f$ is {\it flat} if ${\cal C}(f)\subseteq [-1,1]$.
Note that if $f$ is flat, then
$f$ is also coefficient convex. Typically we denote polynomial coefficients by $c_j$ and $d_j$.\\
\indent The $n$th cyclotomic polynomial $\Phi_n(x)$ (see the next section for
details) has the property that its coefficients tend to be small
in absolute value, e.g., for $n<105$ it is flat. If $n$ has at most three distinct odd prime factors, it can
be shown \cite{GM} that $\Phi_n$ is coefficient convex. 
A question that arises is to what extent the smallness of the coefficients is particular to $\Phi_n(x)$. We
will try to answer this by investigating the coefficients of the other divisors of $x^n-1$ as well. Our
work suggests that as far as the behavior of its coefficients go, $\Phi_n(x)$ does not have a special role
amongst the divisors of $x^n-1$.
Since the number of divisors of $x^n-1$ rapidly increases, we are
only able to say something conclusive in case $n$ has a modest number of divisors. If $n=pq$ or $n=p^2q$,
 then $x^n-1$ has 16, respectively 64 monic divisors (these cases are covered by
 Theorems \ref{plakkerig}, \ref{bluebell}, \ref{natha} and \ref{MAIN}).\\
\indent An exception here is the case where $n$ is a prime power, say $n=p^e$. Then the number of divisors
can get large, but they have a simple structure. Using the uniqueness of the base $p$ representation
Pomerance and Ryan \cite{PR} proved that the divisors of $x^{p^e}-1$ are all flat. We leave it to the
reader to prove the following easy strengthening of this result.
\begin{Thm}
\label{riant}
Let $e\ge 1$ be an integer and $g$ be a monic divisor of $x^{p^e}-1$. We have ${\cal C}(g)=\{1\}$ if $g=(x^{p^j}-1)/(x-1)$ for
some $0\le j\le e$. Furthermore, if $p=2$ and
$g=(x-1)(x^{2^j}-1)/(x^2-1)$, then  for $1\le j\le e$ we have ${\cal C}(g)=\{-1,1\}$. In the remaining cases we have
$${\cal C}(g)=\cases{[0,1] & if $g(1)\ne 0$;\cr
[-1,1] & otherwise.}$$
\end{Thm}
\begin{Thm}
\label{plakkerig}
Let $p<q$ be primes.
Except for $(x-1)\Phi_{pq}(x)$
and $\Phi_p(x)\Phi_q(x)$ all monic divisors of $x^{pq}-1$ are flat. The set of coefficients of $(x-1)\Phi_{pq}(x)$
is of the form $\{-2,-1,1,2\}$ if $p\le 3$ and $[-2,2]$ otherwise.
The set of coefficients of $\Phi_p(x)\Phi_q(x)$ is $[1,\min(p,q)]$.
\end{Thm}
\begin{cor}
All divisors $f\in {\Bbb Z}[x]$ of $x^{pq}-1$ are coefficient convex.
\end{cor}
\begin{Thm} 
\label{bluebell}
Let $p$ and $q$ be distinct primes. Then the monic polynomial divisors of $x^{p^2q}-1$ are 
coefficient convex, with
the exception (in case $q=2$), $(x+1)\Phi_p \Phi_{2p^2}$, where the coefficient set equals $\{-2,0,1,2\}$.
If $\min(p,q)>3$, then all monic divisors - except $x-1$ - are strongly coefficient convex.
\end{Thm}
Let $B(n)$ be the maximum coefficient (in absolute value) that occurs amongst all 
monic divisors of $x^n-1$.
Pomerance and Ryan \cite{PR} conjectured and Kaplan \cite{K2} proved that 
$B(p^2q)=\min(p^2,q)$. Letting $B_{+}(n)$ denote
the maximum amongst all the coefficients of all the monic divisors of $x^n-1$, and $-B_{-}(n)$ the minimum,
we have the following generalization of Kaplan's result.
\begin{Thm}
\label{natha}
Let $p$ and $q$ be distinct primes. Let $1\le p^*\le q-1$
be the inverse of $p$ modulo $q$.
We have $B_{-}(p^2 q)=\min(p,p^*)+\min(p,q-p^*)$ and $B_{+}(p^2 q)=\min(p^2,q)$
\end{Thm}
Note that if $q<p$, then the result gives $B_{-}(p^2q)=B_{+}(p^2q)=q$. (For a more formal definition of $B_{\pm}(n)$ see
Section \ref{height2010}.) The analogue of the latter theorem in case $n=pqr$ is not known, for some
partial results see Kaplan \cite{K2}. Ryan et al. \cite{ryanryan} posed
some conjectures on the basis of extensive numerical calculation.\\
\indent The results stated above (except for Theorem \ref{riant}) are special cases of 
Theorem \ref{MAIN}, our main result, e.g., Theorem \ref{plakkerig} can be read off from
Table 1A. In the derivation of Theorem \ref{natha} we have to use in addition that
$\min(p,p^*)+\min(p,q-p^*)\ge \min(p,q)$. A reformulation of Theorem \ref{MAIN} without tables
is given in Section \ref{compact}.
\begin{Thm}
\label{MAIN}
Let $p$ and $q$ be distinct primes. Let $f(x)\in \Bbb Z[x]$ be a monic divisor of $x^{p^2q}-1$. Then there exists
an integer $0\le k\le 63$ such that 
$$f(x)=f_k(x)=\Phi_1^{k_0}\Phi_p^{k_1}\Phi_q^{k_2}\Phi_{pq}^{k_3}\Phi_{p^2}^{k_4}\Phi_{p^2q}^{k_5},$$
with $0\le k_j\le 1$ and $k=\sum_{j=0}^5k_j2^j$ the binary expansion of $k$. The set of coefficients 
of $f_k$, ${\cal C}(f_k)$, is given in Table {\rm 1}.
\end{Thm}
The difficulty of computing ${\cal C}(f)$ varies rather dramatically; from utterly trivial to challenging in
case of $f_{25}$, $f_{38}$ and $f_{43}$.

\section{Preliminaries}
The $n$th cyclotomic polynomial $\Phi_n(x)$ is defined by
\begin{equation}
\label{gaap}
\Phi_n(x)=\sum_{k=0}^{\phi(n)}a_n(k)x^k=\prod_{d|n}(x^d-1)^{\mu(n/d)},
\end{equation}
where $\mu(n)$ is the M\"obius function and $\varphi(n)$ 
Euler's totient function. Let $p\neq q$ be primes. {}From (\ref{gaap}) we deduce, e.g., that
\begin{equation}
\label{piqu}
\Phi_{pq}(x)={(x-1)(x^{pq}-1)\over (x^p-1)(x^q-1)},
\end{equation}
a formula that will be used repeatedly.\\
\indent We will need the following elementary properties
of $\Phi_n(x)$ (see, e.g., Thangadurai \cite{Thanga} for proofs and a nice introduction
to cyclotomic  polynomials). Throughout we use the letters $p$ and $q$ to denote primes.
\begin{Lem} \label{basic} $~$\\ 
{\rm 1)} $\Phi_n(x)\in \mathbb Z[x]$.\\
{\rm 2)} $\Phi_n(x)$ is irreducible over the rationals.\\
{\rm 3)} $x^n-1=\prod_{d|n}\Phi_d(x)$.\\
{\rm 4)} $\Phi_p(x)=(x^p-1)/(x-1)=1+x+\ldots+x^{p-1}$.\\
{\rm 5)} If $p|n$, then $\Phi_{pn}(x)=\Phi_{n}(x^p)$.\\
{\rm 6)} If $n>1$ is odd, then $\Phi_{2n}(x)=\Phi_n(-x)$.\\
{\rm 7)} For all positive integers $n>1$, we have $\Phi_n(1/x)x^{\phi(n)}=\Phi_n(x)$, that is
$\Phi_n(x)$ is {\tt self-reciprocal}.
\end{Lem}
\noindent For a nonzero polynomial $f\in \mathbb C[x]$, we define its
{\it height} $H(f)$ to be the largest coefficient of $f$ in absolute
value.  For a nonzero polynomial $f\in \mathbb R[x]$, we define $H_{+}(f)$, respectively $H_{-}(f)$ to be the largest, respectively
smallest coefficient of $f$. (In that case $H(f)=\max\{H_{+}(f),|H_{-}(f)|\}$.)
As in \cite{PR}, the observation that if $H(f)=m$, then $H((x^k-1)f(x))\le 2m$
for any positive integer $k$ will be used a few times. We also use that if $f,g\in \mathbb Z[x]$
with deg$(f)\le {\rm deg}(g)$, then 
\begin{equation}
\label{zon}
H(fg)\le (1+{\rm deg}(f))H(f)H(g).
\end{equation}
Another easy observation we need is that if $k>{\rm deg}(f)$, and $m\ge 1$ is an arbitrary
integer, then
\begin{equation}
\label{bier}
{\cal C}_0(f(x)(1+x^k+x^{2k}+\cdots+x^{km}))={\cal C}_0(f).
\end{equation}
If $k>{\rm deg}(f)+1$, then 
${\cal C}(f(x)(1+x^k+x^{2k}+\cdots+x^{km}))={\cal C}(f)\cup \{0\}$. A closely
related observation is that
\begin{equation}
\label{flauuw}
{\cal C}(\Phi_p(x)f(x^p))={\cal C}(f).
\end{equation}
To see this note that if in the coefficient string of $f (=\sum_j c_jx^j)$,
that is in the string $c_0c_1c_2\ldots c_{\rm deg(f)}$, we replace each coefficient by its $p$-fold
repetition (e.g. $c_0c_1$ becomes
$c_0c_0c_0c_1c_1c_1$ if $p=3$), we get the coefficient string of $\Phi_p(x)f(x^p)$.
\subsection{Binary cyclotomic polynomials}
In this subsection we consider the binary cyclotomic polynomials $\Phi_{pq}(x)$ with
$p$ and $q$ distinct primes.\\
\indent In 1883 Migotti proved that  $\Phi_{pq}$ is flat.
Carlitz \cite{CA} noted that if we
drop the zero coefficients in $\Phi_{pq}(x)$, the positive and negative terms occur alternately, 
as, e.g., in
$$\Phi_{21}(x)=x^{12}-x^{11}+x^{9}-x^{8}+x^6-x^4+x^3-x+1.$$
(To prove this,  one can invoke Lemma \ref{tau} below together with (\ref{piqu}).)
Lenstra \cite{Lenstra} (see also Lam and Leung \cite{LL}) gave an explicit description of the coefficients of
$\Phi_{pq}(x)$.
\begin{Lem}
\label{binary} {\rm (\cite{LL}).}
Let $p$ and $q$ be distinct odd primes. Let $\rho$ and $\sigma$ be the (unique) non-negative
integers for which $1+pq=(\rho+1) p+(\sigma+1) q$.
Let $0\le m<pq$. Then either $m=\alpha_1p+\beta_1q$ or $m=\alpha_1p+\beta_1q-pq$
with $0\le \alpha_1\le q-1$ the unique integer such that $\alpha_1 p\equiv m({\rm mod~}q)$
and $0\le \beta_1\le p-1$ the unique integer such that $\beta_1 q\equiv m({\rm mod~}p)$.
The cyclotomic coefficient $a_{pq}(m)$ equals
$$\cases{1 & if $m=\alpha_1p+\beta_1q$ with $0\le \alpha_1\le \rho,~0\le \beta_1\le
\sigma$;\cr -1 & if $m=\alpha_1p+\beta_1q-pq$ with $\rho+1\le \alpha_1\le q-1,~\sigma+1\le 
\beta_1\le p-1$;\cr  0 & otherwise.}
$$
\end{Lem}
The latter lemma does not include the case where $p=2$ and $q$ is odd. However, by
Lemma \ref{basic} we have $\Phi_{2q}(x)=\Phi_q(-x)=1-x+x^2-\cdots+x^{q-1}$.\\
\indent A rather specific observation we will need involving binary cyclotomic polynomials is the 
following.
\begin{Lem}
\label{38-2}
Let $p$ and $q$ be primes with $p<q$. Write $q=fp^2+g$ with $0<g<p^2$ and $\Phi_{pq}(x)=\sum_j a_jx^j$.
We have $$\Phi_{pq}(x)=\sum_{j=0}^{fp-1}a_jx^j+x^{fp}-x^{fp+1}+x^{(f+1)p}+\sum_{j=(f+1)p+1}^{(p-1)(q-1)}a_jx^j.$$
\end{Lem}
{\it Proof}. Write $(p-1)(q-1)=\rho p+\sigma q$ as in Lemma \ref{binary}. Note that $f<\rho$, because otherwise using $\sigma \leq p-2$ and 
$f \geq \rho$ we obtain the contradiction
\begin{eqnarray*}
	(p-1)(q-1)=\rho p+\sigma q & \leq & pq-2q+fp =pq-q-p+((f+1)p-q) \\
	& \leq & pq-p-q<(p-1)(q-1).
	\end{eqnarray*}
Therefore $a_{fp}=1$ and also $a_{(f+1)p}=1$ (Lemma \ref{binary}). Likewise we have $f<q-\rho$, since $f\geq q-\rho$ 
leads to the contradiction
\begin{eqnarray*}
	(p-1)(q-1)=\rho p+\sigma q\geq pq-fp> pq-q/p > pq-q > (p-1)(q-1).
	\end{eqnarray*}
This in combination with the identity $fp+1=(\rho+f+1)p+(\sigma+1)q-pq$ and  Lemma \ref{binary} 
shows that $a_{fp+1}=-1$. Since $(f+1)p<q$, we now see that the terms $x^{fp}-x^{fp+1}+x^{(f+1)p}$ appear consecutively in $\Phi_{pq}(x)$. \qed

\subsection{Inverse cyclotomic polynomials}
\label{tweetwee}
We define $\Psi_n(x)=(x^n-1)/\Phi_n(x)=\sum_{k=0}^{n-\varphi(n)}c_n(k)x^k$ to
be the $n$th {\it inverse cyclotomic polynomial}. It is easy to see, see, e.g., Moree \cite{M},
that $\Psi_1(x)=1$, $\Psi_p(x)=x-1$ and 
\begin{equation}
\label{invi}
\Psi_{pq}(x)=-1-x-x^2-\ldots-x^{\min(p,q)-1}+
x^{\max(p,q)}+x^{\max(p,q)+1}+\ldots+x^{p+q-1}.
\end{equation}
For $n<561$ the polynomials $\Psi_n(x)$ are flat. 
Let $2<p<q<r$ be odd primes.
It is not difficult to show that
$|c_{pqr}(k)|\le [(p-1)(q-1)/r]+1\le p-1$, where $[x]$ denotes the largest {\it integer} $\le x$. Let us call a ternary inverse cyclotomic polynomial
$\Psi_{pqr}(x)$ extremal if for some $k$ we have $|c_{pqr}(k)|=p-1$. Moree \cite{M} showed that 
a ternary inverse cyclotomic polynomial is extremal iff 
$$q\equiv r\equiv \pm 1({\rm mod~}p){\rm ~and~}r<{(p-1)\over (p-2)}(q-1).$$ Moreover, he showed that for an extremal ternary
inverse cyclotomic polynomial $\Psi_{pqr}(x)$ one has ${\cal
C}(\Psi_{pqr})=[-(p-1),p-1]$, and thus that it is strongly coefficient convex. 
\subsection{Inclusion-exclusion polynomials}
Let $\rho=\{r_1,r_2,\ldots,r_s\}$ be a set of natural numbers satisfying $r_i>1$ and $(r_i,r_j)=1$ for $i\ne j$, and put
$$n_0=\prod_i r_i,~n_i={n_0\over r_i},~n_{i,j}={n_0\over r_ir_j}~[i\ne j],\ldots$$
For each such $\rho$ we define a function $Q_{\rho}$ by
$$Q_{\rho}(x)={(x^{n_0}-1)\prod_{i<j}(x^{n_{i,j}}-1)\cdots \over \prod_i (x^{n_i}-1)\prod_{i<j<k}(x^{n_{i,j,k}}-1)\cdots}$$
It turns out that $Q_{\rho}$ is a polynomial, the inclusion-exclusion polynomial. This class of divisors of $x^{n_0}-1$
was introduced by Bachman \cite{B}. He showed that with $D_{\rho}=\{d:d|n_0{\rm~and~}(d,r_i)>1{\rm~for~all~}i\}$, we
have $Q_{\rho}(x)=\prod_{d\in D}\Phi_d(x)$. Furthermore, he showed that ternary $(s=3)$ inclusion-exclusion polynomials
are coefficient convex. Earlier Gallot and Moree \cite{GM} (for alternative proofs, see Bzd{\c e}ga \cite{BZ} and 
Rosset \cite{Rosset}) had shown that in case
$s=3$ and $r_1,r_2,r_3$ are distinct primes, this result is true.
\subsection{On the coefficient convexity of $\Phi_n$ and $\Psi_n$}
\indent In \cite{GM} Theorems \ref{platteband1} and \ref{platteband2} were announced and it was
promised that the present paper would contain the proofs. Here this promise is fulfilled.\\
\indent In \cite{GM} the following result was established. (Its analogue for $\Psi_n$ is
false in general.)
\begin{Thm}
\label{gamo}
{\rm (\cite{GM}).}
Let $n$ be ternary, that is $n=pqr$ with $2<p<q<r$ odd primes. Then, for $k\ge 1$,
$|a_n(k)-a_n(k-1)|\le 1$.
\end{Thm}
It follows that if $n$ is ternary, then $\Phi_n$ is strongly coefficient convex. Using
the latter result one easily proves the following one.
\begin{Thm}
\label{platteband1}
Suppose that $n$ has at most $3$ distinct prime factors, then $\Phi_n$ is 
coefficient convex.
\end{Thm}
{\it Proof}. In case $n$ has at most two distinct odd factors, by Lemma \ref{binary} and Lemma \ref{basic} we
infer that $\Phi_n$ is flat and hence coefficient convex. Now suppose that $n$ is odd. Let $\kappa(n)=\prod_{p|n}p$ be
the square free kernel of $n$. Then, by part 4 of 
Lemma \ref{basic} we have ${\cal C}(\Phi_n)={\cal C}(\Phi_{\gamma(n)})
\cup \{0\}$ if $\kappa(n)<n$. The proof is now completed on invoking Theorem \ref{gamo}. \qed\\

\noindent Numerical computation suggest that if $n$ is ternary, then $\Phi_{2n}$ is coefficient
convex. If this would be true, then in Theorem \ref{platteband1} one can replace `3 distinct prime factors'
by `3 distinct odd prime factors'. This is best possible as the following examples show:\\
$n=7735=5\cdot 7\cdot 13\cdot 17$, $~~{\cal C}(n)=[-7,5]\backslash \{-6\}$\\
$n=530689=17\cdot 19\cdot 31\cdot 53$, $~~{\cal C}(n)=[-50,52]\backslash \{-48,47,48,49,50,51\}$.
\begin{Thm}
\label{platteband2}
Suppose that $n$ has at most $2$, respectively $3$, distinct odd prime factors, then $\Psi_n$ is 
flat, respectively, coefficient convex.
\end{Thm}
{\it Proof}. If $p|n$, then $\Psi_{pn}(x)=\Psi_n(x^p)$. Thus we may restrict to the case where
$n$ is square free. If $n=1$, then $\Psi_1=1$. If $n$ is a prime, then $\Psi_n=x-1$. If
$n$ is composed of two primes, $n=pq$, with $p<q$, then
$$\Psi_{pq}=-1-x-x^2-\cdots-x^{p-1}+x^q+x^{q+1}+\cdots+x^{p+q-1}.$$ 
If $2<p<q$, then
$\Psi_{2pq}=(1-x^{pq})\Psi_{pq}(-x)$. Note that the degree of $\Psi_{pq}$ is smaller than $pq$
and since $\Psi_{pq}(-x)$ is flat, it follows that $\Psi_{2pq}$ is flat.
We conclude that if $n$ has at most two distinct odd prime factors, then $\Psi_n$ is flat. It
remains to consider the case where $n=pqr$, $2<p<q<r$, respectively $n=2pqr$ with $2<p<q<r$.\\
Case 1. $n=pqr$. We have $\Psi_{pqr}(x)=\Phi_{pq}(x)\Psi_{pq}(x^r)$. {}From this identity we
infer that
$$c_{pqr}(k)=\sum_{j=0}^{[k/r]}a_{pq}(k-jr)c_{pq}(j).$$
Put $V_n=\{c_{n}(k):0\le k\le n-\varphi(n)\}$. Choose $k_1$ such that $c_{pqr}(k_1)=\max V_{pqr}=\mu_{+}$.
Then since  $|a_{pq}(k-jr)c_{pq}(j)|\le 1$, we infer that 
$$[1,\mu_{+}]\subseteq \{~c_{pqr}(k_1-jr)~:~0\le j\le [{k_1\over r}]~\}.$$ Similarly one choses $k_2$ such
that $c_{pqr}(k_1)=\min V_{pqr}=\mu_{-}$ and finds that 
$$[\mu_{-},-1]\subseteq \{~c_{pqr}(k_2-jr)~:~0\le j\le [{k_2\over r}]~\}$$ and
hence $V_{pqr}=[\mu_{-},\mu_{+}]$ (by \cite[Lemma 3]{M}
we have $0\in V_{pqr}$). Thus $\Psi_{pqr}$ is coefficient convex.\\
Case 2. $n=2pqr$. A small modification of the above argument gives that $\Psi_{pqr}(-x)$
is coefficient convex. Using that $\Psi_{2n}(x)=(1-x^n)\Psi_n(-x)$ if $n$ is odd and that
$n>n-\varphi(n)={\rm deg}(\Psi_n)$, we infer that also $\Psi_{2pqr}$ is coefficient convex.\\
\indent Thus the proof is completed. \qed\\

\noindent Gallot considered the coefficient convexity of $\Psi_n$  for many $n$ and found that the smallest
$n$ for which it is non-convex is $n=23205=3\cdot 5\cdot 7\cdot 13\cdot 17$. Here the height
is $13$, but $12$ (and $-12$) are not included in ${\cal C}(\Psi_n)$. Further examples (in order of
appearance) are $46410$ (height $13,\pm 12$ not there), $49335$ (height $34$, $\pm 33$ not found),
$50505$ (height $15$, $\pm 14$ not found). There are also examples where a whole range of
values smaller than the height is not in ${\cal C}(\Psi_n)$.

\subsection{Auxiliary polynomials}
In this subsection we determine ${\cal C}(f)$ for various auxiliary polynomials $f$ (where
possible we have adopted the notation of Theorem \ref{MAIN}).
\begin{Lem} 
\label{tau}
Let $u>1$ and $v>1$ be coprime natural numbers.
Put $$\tau_{u,v}(x)={(x-1)(x^{uv}-1)\over (x^u-1)(x^v-1)}.$$
Then $\tau_{u,v}(x)\in \mathbb Z[x]$ is a self-reciprocal flat divisor of $x^{uv}-1$.
If $1<u<v$, then
$${\cal C}(\tau_{u,v})=\cases{\{-1,1\} & if $u=2$;\cr
[-1,1] & otherwise.}$$
The non-negative coefficients of $\tau_{u,v}$ alternate in sign.
\end{Lem}
{\it Proof}. The assumption on $u$ and $v$ ensures that $(x^u-1,x^v-1)=x-1$. Using this
assumption we infer that $\tau_{u,v}(x)\in \mathbb Z[x]$. That $\tau_{u,v}(x)$ is a self-reciprocal divisor
of $x^{uv}-1$ is obvious. We study the coefficients of $\tau_{u,v}(x)$ by first considering the Taylor series
around $x=0$ of the denominator of $\tau_{u,v}(x)$. 
We claim that all coefficients $r_j$ with $j<uv$ in $(1+x^{u}+x^{2u}+\cdots)(1+x^v+x^{2v}+\cdots)=\sum r_jx^j$
are in $[0,1]$. Now if $r_j\ge 2$ and $j<uv$, we can find 
non-negative $\alpha_1,\alpha_2,\beta_1$ and $\beta_2$
such that $j=\alpha_1u+\beta_1v=\alpha_2u+\beta_2v$, with $\alpha_1\ne \alpha_2$ both smaller
than $v$. The latter equality implies however $v|(\alpha_1-\alpha_2)$. This contradiction
completes the proof of the claim. It follows that ${\cal C}(\tau_{u,v})\subseteq [-1,1]$ and
that the non-negative signs alternate. The claim regarding
${\cal C}(\tau_{u,v})$ follows on noting that $\tau_{u,v}=(x^v+1)/(x+1)$ if $u=2$ and
$\tau_{u,v}\equiv 1-x({\rm mod~}x^3)$ if $u\ge 3$. \qed\\

\indent In case $p=3$, the next lemma shows that $\tau_{3,v}(x)$ can be easily given explicitly.
\begin{Lem}
\label{3q}
Let $v$ be a positive integer with $3\nmid v$.
If $v\equiv 1({\rm mod~}3)$, put 
$$f_v(x)=(1-x)(1+x^3+x^6+\cdots+x^{v-1})+x^{v}+(x-1)(x^{v+1}+x^{v+4}+\cdots+x^{2v-3}).$$
If $v\equiv 2({\rm mod~}3)$, put
$$f_v(x)=(1-x)(1+x^3+x^6+\cdots+x^{v-2})+x^{v}+(x-1)(x^{v+2}+x^{v+5}+\cdots+x^{2v-3}).$$
We have $\tau_{3,v}(x)=f_v(x)$.
\end{Lem}
{\it Proof}. Modulo $x^{v}$ we have
$$\tau_{3,v}(x)={(x-1)(x^{3v}-1)\over (x^3-1)(x^v-1)}\equiv (1-x)(1+x^3+x^6+\cdots).$$
We infer that $f_v(x)\equiv \tau_{3,v}(x) ({\rm mod~}x^v)$. We have ${\rm deg}(f_{v})=2v-2={\rm deg}(\tau_{3,v})$, so to finish the proof it
is enough to show that $f_v(x)$ is self-reciprocal (clearly $\tau_{3,v}(x)$ is
self-reciprocal). That is, we have to show that $f_v(1/x)x^{2(v-1)}=f_v(x)$. That this
is the case is easily seen on rewriting $f_v(x)$, in case $v\equiv 1({\rm mod~}3)$ as
$$(1-x)(1+x^3+x^6+\cdots+x^{v-4})+x^{v-1}+(x-1)(x^{v+1}+x^{v+4}+\cdots+x^{2v-3}),$$
and as
$$(1-x)(1+x^3+x^6+\cdots+x^{v-2})+x^{v-1}+(x-1)(x^{v-1}+x^{v+2}+\cdots+x^{2v-3}),$$
in case $v\equiv 2({\rm mod~}3)$. \qed\\

\noindent Lemma \ref{3q} shows that  identical consecutive  coefficients
do not appear in $\tau_{3,v}(x)$ if $(3,v)=1$.
The following lemma determines all polynomials $\tau_{3,v}(x)$ with
this property.
\begin{Lem}
\label{ddrie}
Let $1<u<v$ be coprime integers.
Consecutive coefficients of $\tau_{u,v}(x)$ are always distinct iff $u\le 3$.
\end{Lem}
\begin{cor}
\label{diffidi}
We have $0\in {\cal C}((x-1)\tau_{u,v}(x))$ iff $u > 3$.
\end{cor}
{\it Proof}. If $u=2$ we have $\tau_{2,v}(x)=(x^v+1)/(x+1)$ and so consecutive coefficients
are always distinct. If $u=3$ it is seen from Lemma \ref{3q} that this property
also holds. Proceeding as in the proof of Lemma \ref{3q} we find that modulo $x^v$ we
have $\tau_{u,v}(x)\equiv (1-x)(1+x^u+x^{2u}+\cdots)$ and hence, if $u\ge 4$, the second and third
coefficient of $\tau_{u,v}(x)$ both equal zero. \qed

\begin{Lem}
\label{nulgat}
Let $1<u<v$ be coprime numbers. Put $h=(x-1)\tau_{u,v}(x)$.
We have 
$${\cal C}(h)=\cases{\{-2,-1,1,2\} & if $u\le 3$;\cr
[-2,2] & otherwise.}$$
\end{Lem}
{\it Proof}. Put $d=(u-1)(v-1)$. Using the self-reciprocity of $\tau_{u,v}(x)$  we infer that
$\tau_{u,v}(x)=x^d-x^{d-1}+\cdots-x+1$. On writing
$h(x)=\sum_j c_jx^j$, we now deduce that 
$c_0=-1$, $c_1=2$, $c_{d}=-2$ and $c_{d+1}=1$. Since clearly
${\cal C}(h)\subseteq [-2,2]$ (use Lemma 
\ref{tau}), we infer that $\{-2,-1,1,2\}\subseteq {\cal C}(h)\subseteq [-2,2]$.
On invoking Corollary \ref{diffidi}, the proof is then completed. \qed

\begin{Lem}
\label{vijf}
Let $u,v$ be natural numbers.
Put $$\sigma_{u,v}(x)={(x^u-1)\over (x-1)}{(x^v-1)\over (x-1)}=\sum_{j=0}^{u+v-2}c_jx^j$$ 
W.l.o.g. assume that $u\le v$. We have
$$
c_j=\cases{
j+1 & if $0\le j\le u-1$;\cr
u & if $u\le j\le v-1$;\cr
v+u-j-1 & if $v\le j\le v+u-2$.}
$$
It follows that
${\cal C}(\sigma_{u,v})=[1,u]$. If $(u,v)=1$, then
$\sigma_{u,v}(x)|x^{uv}-1$.
\end{Lem}
\begin{cor}
\label{ojoj}
If $u<v$, then ${\cal C}((x-1)\sigma_{u,v}(x))=[-1,1]$.
\end{cor}
\begin{cor}
If $(u,v)=1$, then $B(uv)\ge B_{+}(uv)=\min(u,v)$.
\end{cor}
\begin{cor}
\label{22}
Put $f_{22}=\Phi_p\Phi_q\Phi_{p^2}$. Then ${\cal C}(f_{22})=[1,\min(p^2,q)]$.
\end{cor}
{\it Proof of Lemma} \ref{vijf}. Modulo $x^{u}$ we have
$$\sigma_{u,v}(x)\equiv {1\over (1-x)^2}\equiv \sum_{j=1}^{u}jx^{j-1} ({\rm mod~}x^{u}),$$
showing that $c_j=j+1$ if $0\le j\le u-1$. That $c_j=u$ if $u\le j\le v-1$ is obvious. Using
that $\sigma_{u,v}$ is self-reciprocal, it then follows that $c_j=v+u-j-1$ if $v\le j\le v+u-2$.\\
\indent If $(u,v)=1$, then $((x^u-1)/(x-1),(x^v-1)/(x-1))=1$ and using this we
infer that $\sigma_{u,v}(x)|x^{uv}-1$. \qed

\begin{Lem}
\label{20}
Let $p$ and $q$ be distinct primes. Put $f_{20}=\Phi_q\Phi_{p^2}$ We have
$${\cal C}(f_{20})=\cases{[1,\min([{q-1\over p}]+1,p)] & if $p< q$;\cr
[0,1] & if $p>q$.}$$
Consequently $f_{20}$ is flat iff $p>q$.
\end{Lem}
{\it Proof}. Left as an exercise to the interested reader.\qed

\begin{Lem}
\label{abc}
Let $a,b,c$ be positive integers. Put
$$g_{a,b,c}(x)=(1+x+\ldots+x^{a-1}+2x^{a}+\ldots+2x^{a+b-1})(1+x+x^2+\ldots+x^{c-1}).$$
Alternatively one can write
$$g_{a,b,c}(x)=\Big({2x^{a+b}-x^a-1\over x-1}\Big)\Big({x^c-1\over x-1}\Big).$$
Suppose $a$ is odd. Then $g(=g_{a,b,c})$ is coefficient convex. We  have
${\cal C}(g)=[1,\mu]$, with
$$\mu=\cases{2c & if $c\le b$;\cr
\min(b+c,a+2b) & if $c>b$.}$$ 
\end{Lem}
\begin{cor}
\label{bogdan}
Put ${\overline g}=x^{a+b+c-2}g_{a,b,c}(1/x)$. We have
$${\overline g}={\overline g}_{a,b,c}=\Big({x^{a+b}+x^b-2\over x-1}\Big)\Big({x^c-1\over x-1}\Big).$$
If $a$ is odd, then $\overline{g}$ is coefficient convex and ${\cal C}(\overline{g})=[1,\mu]$.
\end{cor}
{\it Proof of Lemma} \ref{abc}. To find the maximum coefficient of $g$ is easy. It is the coefficient
convexity that is slightly less trivial. Write $g=\sum_{j=0}^{a+b+c-2}d_jx^j$. 
We consider two cases.\\
Case 1. $c\ge a+b$. We have to show that all coefficients $1,2,\ldots,\mu$, where
$\mu=a+2b$, occur. It
is easy to see that $\{d_0,\ldots,d_{a+b-1}\}$ contains all odd number $\le \mu$ (here
we use the assumption that $a$ is odd). Likewise one sees that 
$\{d_{c},\ldots,d_{a+b+c-2}\}$ contains all even integers $\le \mu$.\\
Case 2. $c<a+b$. Here we proceed by induction with respect to $c$. For $c=1$ we have
1 and 2 as coefficients and we are done. Suppose the result is true up to $c_1$. We
want to show it for $c=c_1+1$. Here at most two new coefficient values can arise, namely
the previous maximum, $\mu_{c_1}$, with 1 added and the previous maximum with 2 added. In the latter
case (which only arises if $c\le b$) we have to show that $\mu_{c_1}+1$ also occurs as coefficient. 
The coefficient of $d_{a+c-1}=2c$ is the new maximum here. Note that $d_{a+c-2}=2c-1$.
Thus using the induction hypothesis the set of coefficients equals $\{1,2,\ldots,\mu_{c_1},
\mu_{c_1}+1,\mu_{c_1}+2\}$ and is hence coefficient convex. \qed

By $[f]_{x^k}$ we denote the coefficient of $x^k$ in $f$.
\begin{Lem}
\label{koffie}
Let $p$ and $q$ be distinct primes. Put $f_{24}=\Phi_{pq}\Phi_{p^2}$.
Let $1\le p^*\le q-1$
be the inverse of $p$ modulo $q$. Write $f_{24}=\sum_k c_kx^k$.\\
{\rm 1)} We have 
$${\cal C}(f_{24})=\cases{[-\min(q-p^*,p),\min(p^*,p)] & if both $p$ and $q$ are odd;\cr
[-\min(q-p^*,p),\min(p^*,p)]\backslash \{0\} & otherwise.}$$
Consequently, $f_{24}$ is flat iff $q=2$.\\
{\rm 2)} Let $k\ge 0$ and $\min(p,q)>2$. We have $c_{1+kp}=-[\sigma_{q-p^*,p}(x)]_{x^k}$ and $c_{kp}=[\sigma_{p^*,p}(x)]_{x^k}$. 
If $2p^*<q$, then $c_{2+kp}=[x^{q-2p^*}\sigma_{p^*,p}(x)]_{x^k}$. If $2p^*>q$, then
$c_{-1+kp}=-[x^{2p^*-q}\sigma_{q-p^*,p}(x)]_{x^k}$.
\end{Lem}
{\it Proof}. 1) The case where $p$ or $q$ is even is left to the reader. So let us assume
that both $p$ and $q$ are odd.
The $k$th coefficient $c_k$ in $f_{24}$ equals 
$$\sum_{k\ge 0\atop 0\le k-jp<pq,~0\le j\le p-1}a_{pq}(k-jp).$$
Since this is a sum of binary cyclotomic coefficients by Lemma \ref{binary} we
have 
$$-(q-1-\rho)\le c_k\le \rho+1{\rm ~and~}-p\le c_k\le p.$$ 
On noting that $\rho+1=p^*$ we thus
obtain that $-m_2\le c_j\le m_1$ with 
$m_2=\min(q-p^*,p)$ and $m_1=\min(p^*,p)$. Using Lemma \ref{binary}  
we obtain that $c_{jp}=\sum_{j_1=0}^j a_{pq}(j_1p)=j+1$ for 
$0\le j\le m_1-1$. Likewise we find on using that $1=(\rho+1)p+(\sigma+1)q-pq$ that
$c_{jp+1}=-j-1$ for $0\le j\le m_2-1$. Since $f_{24}\equiv 1-x ({\rm mod~}x^3)$, it
follows that $0\in {\cal C}(f_{24})$.\\
2) Note that $c_{1+kp}$ is the coefficient of $x^{1+kp}$ in
$$\Phi_p(x^p)\sum_{0\le j<q}a_{pq}(1+jp)x^{1+jp}.$$
Using Lemma \ref{binary} we then infer that the latter polynomial
equals
$$-x\Big({x^{p(q-p^*)}-1\over x^p-1}\Big)\Big({x^{p^2}-1\over x^p-1}\Big).$$
It follows that $c_{1+kp}$ is the coefficient of $x^k$ in $-\sigma_{q-p^*,p}(x)$. 
A similar argument gives $c_{kp}=[\sigma_{p^*,p}(x)]_{x^k}$.
{}From $1+pq=p^*p+q^*q$ we obtain $2=2p^*p+(2q^*-p)q-pq$. The assumption $2p^*<q$
implies $q^*>p/2$ and hence $1\le 2p^*<q$ and $1\le 2q^*-p<q^*$.
Reasoning as before we then find that $c_{2+kp}$ is the coefficient of $x^k$ in $x^{q-2p^*}\sigma_{p^*,p}(x)$.
Likewise the final assertion is established. \qed

\begin{Lem} 
\label{blue}
Put $f_{25}=(x-1)\Phi_{pq}\Phi_{p^2}$. Define $\gamma(p,q)=\min(p,p^*)+\min(p,q-p^*)$. 
Suppose $\min(p,q)>2$.
Write
$f_{25}=\sum d_jx^j$. We have $\{d_{1+kp}\}_{k=0}^{\infty}=[0,\gamma(p,q)]$.
If $2p^*<q$, then $\{d_{2+kp}\}_{k=0}^{\infty}=[-\gamma(p,q),0]$. 
If $2p^*>q$, then $\{d_{kp}\}_{k=0}^{\infty}=[-\gamma(p,q),0]$.
\end{Lem}
{\it Proof}. Using part 2 of Lemma \ref{koffie} we find that $d_{1+kp}=c_{kp}-c_{1+kp}=[\sigma_{p^*,p}(x)+\sigma_{q-p^*,p}(x)]_{x^k}$. Note that
$$\tau(x):=\sigma_{p^*,p}(x)+\sigma_{q-p^*,p}(x)=\Big({x^{q-p^*}+x^{p^*}-2\over x-1}\Big)\Big({x^p-1\over x-1}\Big).$$
We have
$$\tau(x)=\cases{{\overline g}_{q-2p^*,p^*,p}(x) & if $q>2p^*$;\cr
{\overline g}_{2p^*-q,q-p^*,p}(x) & if $q<2p^*$.}$$
On invoking 
Corollary \ref{bogdan} we then obtain, after an easy computation to verify that $\mu=\gamma(p,q)$, that
$\{d_{1+kp}\}_{k=0}^{\infty}={\cal C}(\tau)\cup \{0\}=[0,\gamma(p,q)]$.\\
\indent Using part 2 of Lemma \ref{koffie} and the assumption $q>2p^*$, we find that
$$d_{2+kp}=-c_{2+kp}+c_{1+kp}=-[x^{q-2p^*}\sigma_{p^*,p}(x)+\sigma_{q-p^*,p}(x)]_{x^k}.$$
Now
$$x^{q-2p^*}\sigma_{p^*,p}(x)+\sigma_{q-p^*,p}(x)=\Big({2x^{q-p^*}-x^{q-2p^*}-1\over x-1}\Big)\Big({x^p-1\over x-1}\Big)
=g_{q-2p^*,p^*,p}(x).$$
Using Lemma \ref{abc} we obtain
that $\{d_{2+kp}\}_{k=0}^{\infty}={\cal C}(-\tau)\cup \{0\}=[-\gamma(p,q),0]$.\\
\indent The proof of the final assertion is similar and left to the reader. \qed

\begin{Lem}
\label{hulpp}
Let $q>3$ be a prime. Then the coefficients of the polynomial \\ $g:=(x-1)(1+x^3+x^6)\Phi_{3q}(x)$ are all nonzero.
\end{Lem}
{\it Proof}. Since $g=\sum c_jx^j$ is anti-self-reciprocal, it suffices to show
that $c_j\ne 0$ for $0\le j\le q+2$. Modulo $x^{2q}$, we have
$$g\equiv -(1-2x+x^2)(1+2x^3+3\sum_{j=2}^{\infty}x^{3j})(1+x^q),$$ and so clearly
$c_0,c_1,\ldots,c_{q-1}$ are all nonzero. By computation one checks that also $c_q$, $c_{q+1}$ and
$c_{q+2}$ are nonzero. Alternatively the proof is completed on noting that the sum of any two
coefficients in $(1-2x+x^2)(1+2x^3+3\sum_{j=2}^{\infty}x^{3j})$ that are $q$ apart (here we use
that $q\ge 5$) is nonzero. \qed

\begin{Lem}
\label{hulpq}
Let $p>3$ be a prime. Then $0\in {\cal C}((x-1)\Phi_{3p}\Phi_{p^2})$.
\end{Lem}
{\it Proof}. Put $f(x)=(x-1)\Phi_{3p}\Phi_{p^2}$.
If $p\equiv 1({\rm mod~}3)$, then by Lemma \ref{3q} we find that
$$f(x)\equiv -(1-x)^2(1+x^3+\ldots+x^{p-4})-x^{p-1}({\rm mod~}x^{p+1}),$$
and hence $c_p=0$. If $p\equiv 2({\rm mod~}3)$, then by Lemma \ref{3q} we find that
$$f(x)\equiv -(1-x)^2(1+x^3+\ldots+x^{p-2})-2x^p+3x^{p+1}({\rm mod~}x^{p+3}),$$
and hence $c_{p+2}=0$. \qed

\begin{Lem}
\label{25}
Put $f_{25}=(x-1)\Phi_{pq}\Phi_{p^2}$. Define $\gamma(p,q)=\min(p,p^*)+\min(p,q-p^*)$. Then
$${\cal C}(f_{25})=\cases{[-\gamma(p,q),\gamma(p,q)]\backslash \{0\} & if $p\le 3$ and $q\neq 2$;\cr
[-\gamma(p,q),\gamma(p,q)] & otherwise.}$$
Consequently, $f_{25}$ is never flat.
\end{Lem}
{\it Proof}. Note that if ${\cal C}(f)\subseteq [-a,b]$ with $a$ and $b$ non-negative, then 
${\cal C}((x-1)f)\subseteq [-a-b,a+b]$.
By Lemma \ref{koffie} we thus infer that ${\cal C}(f_{25})\subseteq [-\gamma(p,q),\gamma(p,q)]$.\\
If $q=2$, then $\gamma(p,2)=2$ and one easily sees that ${\cal C}(f_{25})=[-2,2]$.\\
If $p=2$ and $q=3$, then ${\cal C}(f_{25})=[-\gamma(2,3),\gamma(2,3)]\backslash \{0\}=[-3,3]\backslash \{0\}$.\\
If $p=2$ and $q>3$, then the coefficients of $f_{24}$ are alternating in sign and so $0\not\in {\cal C}(f_{25})$. 
The consecutive coefficients are $-1,2,-3,-4,\ldots,-4,3,-2,1$, where the coefficients not indicated are all
$\pm 4$. One
infers that  $${\cal C}(f_{25})=[-\gamma(2,q),\gamma(2,q)]\backslash \{0\}=[-4,4]\backslash \{0\}.$$
So we have dealt with the 
case $\min(p,q)=2$ and may assume that $\min(p,q)>2$. 
Then from  ${\cal C}(f_{25})\subseteq [-\gamma(p,q),\gamma(p,q)]$ and
Lemma \ref{blue} we conclude that ${\cal C}_0(f_{25})=[-\gamma(p,q),\gamma(p,q)]$. It remains to determine whether
$0\in {\cal C}(f_{25})$. \\
If $\min(p,q)>3$, then the coefficient of $x^3$ is zero, so assume that $\min(p,q)=3$.\\
If $p=3$, then by Lemma \ref{hulpp} we see that $0\not\in {\cal C}(f_{25})$.\\ 
If $q=3$, then by 
Lemma \ref{hulpq} we see that $0\in {\cal C}(f_{25})$.
\qed

\begin{Lem}
\label{26}
Let $p$ and $q$ be distinct primes. Put $f_{26}=\Phi_p\Phi_{pq}\Phi_{p^2}$ and $f_{27}=(x-1)f_{26}$.
Then ${\cal C}(f_{26})=[0,1]$ and ${\cal C}(f_{27})=[-1,1]$.
\end{Lem}
{\it Proof}. Write $f_{26}=\sum_j c_jx^j$ and $f_{27}=\sum_j d_jx^j$. 
Note that $f_{26}=(\Phi_p\Phi_{pq})\Phi_p(x^p)=\Phi_p(x^q)\Phi_p(x^p)$ and thus $f_{26}$ has only non-negative
coefficients. Since the equation $aq+bp=a'q+b'p$ with $a,a'\le p-1$ has only the solution $a=a'$ and $b=b'$ it
follows that ${\cal C}(f_{26})\subseteq [0,1]$. On checking that $c_0=1$ and $c_1=0$ it follows that 
${\cal C}(f_{26})=[0,1]$ and hence ${\cal C}(f_{27})\subseteq [-1,1]$. Note that
$d_0=-1$, $d_1=1$. Using that, in case $q=2$, 
$$-f_{27}\equiv {x^p+1\over x+1} ({\rm mod~}x^{p+1}),$$
we easily compute that $d_j=0$ with
$$j=\cases{3 & if $p=2$, $q=3$;\cr
4 & if $p=2$, $q>3$;\cr
p & $q=2$, $p\ge 3$;\cr
2 & if $p\ge 3$, $q\ge 3$.}$$ 
This concludes the proof. \qed

\begin{Lem}
\label{30}
Let $p$ and $q$ be distinct primes. Put $f_{30}=\Phi_p\Phi_q\Phi_{pq}\Phi_{p^2}$.
We have $${\cal C}(f_{30})=[1,\min(p,q)].$$
\end{Lem}
{\it Proof}. Note that $f_{30}=(1+x+\ldots+x^{pq-1})(1+x^p+\ldots+x^{(p-1)p})$. Write
$f_{30}=\sum c_kx^k$. We have
$$0\le c_k=\sum_{0\le k-jp<pq\atop 0\le j\le p-1}1\le \min(p,q).$$
For $0\le r\le \min(p,q)-1$ we have $c_{rp}=r+1$. It is easy to see that
$0$ is not in ${\cal C}(f_{30})$.\qed

\begin{Lem}
\label{f36}
We have ${\cal C}(f_{36})=[-1,1]$. 
\end{Lem}
{\it Proof}. Rewriting shows that $f_{36}(x)=\Phi_q(x)\Phi_{pq}(x^p)$. Because of the alternating 
character of the coefficients of $\Phi_{pq}$ after dropping the zeros, we immediately conclude that $H_+(f_{36})=1$ and $H_-(f_{36})\geq-1$. It is 
also obvious that we have $H_-(f_{36}) = -1$ if $p>q$. In case $p<q$ we express $f_{36}$ 
differently:
\begin{eqnarray}
f_{36}(x) & = & \Phi_q(x)\Phi_{pq}(x^p)=\frac{(x^q-1)}{(x-1)}\cdot
\frac{(x^{p^2q}-1)}{(x^{pq}-1)}\cdot\frac{(x^p-1)}{(x^{p^2}-1)}.\nonumber
\end{eqnarray}
Using the power series for $(1-x^{p^2})^{-1}$ we obtain
\begin{eqnarray}
\label{bloblo}
f_{36}(x) & = &\frac{(x^p-1)(x^q-1)}{1-x}\cdot\frac{x^{p^2q}-1}{x^{pq}-1}\cdot\frac 1{1-x^{p^2}}\nonumber\\
 & = & (1+x+\ldots+x^{p-1}-x^q-x^{q+1}-\ldots-x^{p+q-1})\cdot\nonumber\\ &   & (1+x^{pq}+
\ldots+x^{(p-1)pq})\cdot(1+x^{p^2}+x^{2p^2}+\ldots).
\end{eqnarray}
Let us assume that $H_-(f_{36})>-1$.\\
Since $pq>p+q$ the second factor in (\ref{bloblo}) can be neglected modulo $x^{p+q}$. Hence we have 
$$[f_{36}(x)]_{x^q} = -1 + [(1+x+\ldots+x^{p-1})(1+x^{p^2}+x^{2p^2}+\ldots)]_{x^q}.$$
Now our assumption implies that $q$ can be written as $n \cdot p^2 + r$ with $1 \leq r \leq p-1$, because otherwise $[f_{36}(x)]_{x^q} = -1$. But then we have
$$[f_{36}(x)]_{x^{p+q-1}} = -1 + [(1+x+\ldots+x^{p-1})(1+x^{p^2}+x^{2p^2}+\ldots)]_{x^{p+q-1}}.$$
Our assumption implies $p+q-1 = n' \cdot p^2 + r'$ with $0 \leq r' \leq p - 1$.
With $p^2>(p-1)+(p-1)$, we conclude $n=n'$ and hence $p-1=r'-r$. But we have $r'-r \leq (p - 1) - 1$.
Therefore the assumption that $H_-(f_{36})>-1$ must be false and we conclude $H_-(f_{36})=-1$.\\
\indent {}From (\ref{bloblo}) we infer that the coefficient of $x^p$ is zero if $p<q$. If $p>q$, then
clearly $f_{36}\equiv \Phi_q(x) ({\rm mod~}x^{q+1})$ and the coefficient of $x^q$ is zero. We conclude
that the coefficient of $x^{\min(p,q)}$ is zero and hence the proof is completed.\qed\\

The next three lemmas will be used in order to establish Lemma \ref{38}.
\begin{Lem}
\label{38flauw}
Let $p$ and $q$ be distinct primes. Put $f_{38}=\Phi_p\Phi_{q}\Phi_{p^2q}$. We have
${\cal C}(f_{38})\subseteq [-\min(p,q),\min(p,q)]$.
\end{Lem}
{\it Proof}. Note that $f_{38}=\Phi_pf_{36}=\Phi_qf_{34}$. On using that $H(f_{34})=1$ (easy on using (\ref{flauuw}))
and  $H(f_{36})=1$ (by Lemma \ref{f36}) and invoking (\ref{zon}), it follows that $H(f_{38})\le \min(p,q)$. \qed

\begin{Lem}
\label{38-1}
Let $p$ and $q$ be distinct odd primes. Put $f_{38}=\Phi_p\Phi_{q}\Phi_{p^2q}$ and $\beta(p,q)=\min(p,q,q ({\rm mod~}p^2),p^2-q ({\rm mod~}p^2))$.
We have $[-\beta(p,q),0]\subseteq {\cal C}(f_{38})$.
\end{Lem}
{\it Proof}. Write $f_{38}=\sum_j d_jx^j$.\\
1) The case $q<p$. Here we have $\beta(p,q)=q$.\\
Observe that
\begin{eqnarray*}
	f_{38}(x) & = & (1+2x+\ldots+qx^{q-1}+\ldots+qx^{p-1}+\ldots+x^{p+q-2})(1-x^p+x^{pq}-\ldots)\\
	& = & 1+2x+\ldots+qx^{q-1}+\ldots+qx^{p-1}+\ldots+x^{p+q-2}\\
	&   & -x^p-2x^{p+1}-\ldots-qx^{p+q-1}-\ldots-qx^{2p-1}-\ldots-x^{2p+q-2}+x^{pq}\pm\ldots
\end{eqnarray*}
Furthermore, since $p>q \geq 3$ we have $2p+q-1<3p-1<pq$, and hence the coefficients $[-q,0]$  appear 
from $-qx^{2p-1}$ to $0x^{2p+q-1}$.\\
2) The case $p<q$.\\
We have
$$
f_{38}(x) = (1+2x+\ldots+px^{p-1}+\ldots+px^{q-1}+\ldots+x^{p+q-2})(1-x^p+x^{p^2}-\ldots).
$$
Let $0\leq y \leq \beta(p,q)$. We need to show that $-y \in {\cal C}(f_{38})$.\\
Since $f_{38}(x)=(1+2x+\ldots+px^{p-1}+\ldots+px^{q-1}+\ldots+x^{p+q-2})\Phi_{pq}(x^p)$, we have to evaluate the combinations of these two factors to get the coefficients. So we may express $d_j$, the $j$th coefficient of 
$f_{38}$, as
$$d_j=\sum_{k\ge 0 \atop 0\leq j-kp \leq p+q-2}{\min\{j-kp+1,p,p+q-1-j+kp\}\cdot a_k}.$$
Let
$$
e_{j,k}=\cases{\min\{j-kp+1,p,p+q-1-j+kp\} & if $0 \leq j-kp \leq p+q-2$; \cr
 0 & otherwise,}
$$
so 
\begin{equation}
\label{deee}
d_j=\sum_{k\ge 0}{a_k e_{j,k}}.
\end{equation}
We make the following observation:\\
{\tt Observation}: If $e_{j,k}\neq 0$ and $e_{j,k+l} \neq 0$, we have $e_{j,k+i}=p$ for all $0<i<l$.

i) Assume $y \geq 2p-g$ and $y\neq p$.\\
%Since $p$ and $q$ are prime, we know $g\neq p$. If $g<p$, then $\beta(p,q)=g<p$ and hence $y \neq p$.
%If $g>p$, we have $y \geq 2p-g < p$, so regardless of $g$, we have $y \neq p$.\\
We want to determine $d_{y+q-1}$. Since $y\leq p^2-g$, we have $$y+q-1<p^2-g+fp^2+g=(f+1)p^2,$$
and hence $e_{y+q-1,kp}=0$ for $k\ge (f+1)p$.
By Lemma \ref{38-2} $x^{fp}-x^{fp+1}+x^{(f+1)p}$ appear consecutively 
in $\Phi_{pq}(x)$, i.e. $a_i=0$ for $fp+1<i<(f+1)p$, the largest index $k$ for which $a_k e_{y+q-1,k}$ is not vanishing is $k=fp+1$. We have 
$$e_{y+q-1,fp+1} = \min\{q+y-fp^2-p,p,p-y+fp^2+p\} = \min\{y+g-p,p\} = p,$$
since $y+g-p\geq 2p-g+g-p=p$ by the assumption on $g$. Furthermore, 
$$e_{y+q-1,0} = \min\{y+q,p,p-y\} = p-y > 0.$$
Using the Observation and the alternating character of the non-negative coefficients $a_k$, we conclude that
\begin{eqnarray*}
d_{y+q-1} & = & a_0 \cdot (p-y) + \sum_{k=1}^{fp}{a_k e_{y+q-1,k}}+ a_{fp+1} \cdot p\\
& = & (p-y) - p + p - \ldots + p - p = - y.
	\end{eqnarray*}

ii) Assume $y \geq 2p-(p^2-g)$ or $y=p$.\\
If $y=p$, then $\beta(p,q)=p$ and hence $p^2-g>p$, so $y = p \geq 2p-(p^2-g)$.\\
This time we want to calculate $d_{2p+q-y-1}$. Because $2p+q-y-1>p+q-2$, we have $e_{2p+q-y-1,0}=0$. However, we have $a_{1}=-1$ and $$e_{2p+q-y-1,1}= \min\{p+q-y,p,y\}=y.$$ Since $y\geq 2p-(p^2-g)$, we have 
$$2p+q-y-1<2p+q-2p+(p^2-g)\leq fp^2+g+(p^2-g)=(f+1)p^2,$$ so $e_{2p+q-y-1,(f+1)p}=0$. Furthermore,
\begin{eqnarray*}
e_{2p+q-y-1,fp+1} & = & \min\{2p+q-y-fp^2-p,p,y-p+fp^2+p\}\\
& = & \min\{p+g-y,p,y+fp^2\} = \min\{p,y+fp^2\}.
	\end{eqnarray*}

iia) Consider the case that $y>0$.\\
The smallest index $k$ for which $a_k e_{2p+q-y-1,k}\ne 0$  is $1$, since $e_{2p+q-y-1,1}=y>0$. Using the considerations above, the largest index is $k=fp+1$ (use Lemma \ref{38-2}).
If $f=0$, we have $1=fp+1$, so we directly conclude $d_{2p+q-y-1} = a_1\cdot y = - y$.
If $f > 0$, we have $fp^2>p$, so we have 
%$$e_{2p+q-y-1,fp+1} = \min\{2p+q-y-fp^2-p,p,y-p+fp^2+p\} = \min\{p+g-y,p,y+fp^2\}$$
%and hence
$e_{2p+q-y-1,fp+1} = p$. Reasoning as under i) we obtain
\begin{eqnarray*}
d_{2p+q-y-1} & = & a_1\cdot y + \sum_{k=2}^{fp}{a_k e_{2p+q-y-1,k}}+ a_{fp+1} \cdot p\\
& = & - y + p - p + \ldots + p - p = - y.
	\end{eqnarray*}

iib) Now assume that $y=0$.\\
We have $e_{2p+q-y-1,0}=0$ and $e_{2p+q-y-1,1}= \min\{p+q-y,p,y\}=y=0$. We know that $a_k=0$ for $1<k<p$. Now let $k\geq p$.\\
\indent If $f=0$, i.e. $q<p^2$, we have $2p+q-y-1<(f+1)p^2=p^2$ and hence $e_{2p+q-y-1,k}=0$, since $(2p+q-y-1)-kp < 0$ . So we have established, that the sum $\sum_{k}{a_k e_{2p+q-y-1,k}}$ only consists of summands which are zero, hence $d_{2p+q-y-1}=0=y$.\\
\indent Let $f>0$. The smallest index $k$ for which $a_k e_{2p+q-y-1,k}$ does not vanish is $p$. For $k<p$ we 
have seen before that $a_k e_{2p+q-y-1,k}$ vanishes and further
$a_p=1$ as well as $e_{2p+q-y-1,p}=\min\{2p+q-p^2,p,p^2-p\}=p$.
The largest index $k$ for which $a_k e_{2p+q-y-1,k}$ does not vanish is $fp+1$, with the same reasoning as 
in iia) for $f>0$ and again $e_{2p+q-y-1,fp+1}=p$. So we obtain
\begin{eqnarray*}
d_{2p+q-y-1} & = & a_p\cdot p + \sum_{k=p+1}^{fp}{a_k e_{2p+q-y-1,k}}+ a_{fp+1} \cdot p\\
& = & p - p + \ldots + p - p = 0 = - y.
	\end{eqnarray*}

iii) Consider $y$ to be arbitrary.\\
If $y$ does not have the properties of case i) or case ii), then $0\leq 2y \leq 4p - p^2 - 2$ by a simple addition. It is easy to see that $p=3$ and $y=0$ is the only possibility (with $p$ odd) for that to happen.
We will show that if $p=3$, then  $0\in {\cal C}_{\le 0}(f_{38})$. We have
$$
f_{38}(x) = \Phi_3(x)\Phi_q(x)\Phi_{3q}(x^3)= (1+2x+3x^2\ldots+3x^{q-1}+2x^q+x^{q+1})(1-x^3+x^9\ldots)
$$
If $q=5$ then $d_9=0$ by direct computation. If $q\geq 7$ it is easy to see that $d_5=0$.\\
\indent If $0\le y\le \beta(p,q)$, then we are either in case i), case ii) or $p=3$ and $y=0$. All these
cases have been dealt with. \qed

\begin{Lem}
\label{38-0}
Let $p$ and $q$ be distinct odd primes. 
We have $\min {\cal C}(f_{38})\ge -\beta(p,q)$.
\end{Lem}
{\it Proof}. Note that $\beta(p,q)\le \min(p,q)$. If $\beta(p,q)=\min(p,q)$ we are done by Lemma \ref{38flauw}, so assume that $\beta(p,q)<\min(p,q)$.
Write $q=fp^2+g$, $0<g<p^2$. Note that $\beta(p,q)<\min(p,q)$ implies that $p<q$ and 
$0<g<p$ or $p^2-p<g<p^2$ and so we have to show that all other coefficients are larger 
than or equal to $-g$ and $-(p^2-g)$, respectively.\\
We will use
$$
\Phi_{pq}(x)  =  \frac{(x^{pq}-1)(x-1)}{(x^p-1)(x^q-1)} = (1-x)(1+x^p+\ldots+x^{p(q-1)})(1+x^q+x^{2q}\ldots).
$$
Let $S$ be the numerical semigroup generated by the primes $p$ and $q$, that is the set of all linear combination of $p$ and $q$ of the form $mp+nq$ with $m,n\ge 0$.\\
Note that (cf. \cite{GM})
$$a_k=\cases{1 & if $k\in S$ and $k-1\not\in S$;\cr
-1 & if $k\not\in S$ and $k-1\in S$\cr 
0 & otherwise.}$$
We write $f_{38}(x)=\sum_j d_jx^j$ and use the notation $e_{j,k}$ that was introduced in Lemma \ref{38-1}.
To get a lower bound for an arbitrary coefficient $d_j$, let $n$ be the smallest number such that $a_n=-1$ and $0\leq j-np\leq p+q-2$, i.e. such that $a_n e_{j,n}<0$.
Likewise we let $N$ be the largest number such 
that $0\leq j-Np\leq p+q-2$ and $a_N=-1$. If $n$ and $N$ do not exist, then 
$d_j\geq 0$ and we are done, so now assume that $n$ and $N$ exist and put $d=N-n$.\\
If $d>0$ (and hence $n<N$), we have $e_{j,n}\leq p+q-1-j+np$ and $e_{j,N}\leq j-Np+1$ and thus, using (\ref{deee}),
\begin{eqnarray*}
d_j & = & \sum_{k}{a_k e_{j,k}} \geq \sum_{k=n}^N{a_k e_{j,k}}\\
 & = & -e_{j,n} + p - p \pm \ldots + p -e_{j,N}\\
 & \geq & (-p-q+1+j-np) + p + (Np-j-1) = Np-np-q.
\end{eqnarray*}
If $d=0$, we have $d_j\geq -q =dp-q$, so we infer that always
\begin{equation}
\label{-g}
d_j\geq -fp^2-g+dp.
\end{equation}
The above inequality does not suffice to deal with small $d$. To this end we will need the following two claims.\\
{\tt Claim 1}: Let $m\ge 1$ be an arbitrary integer. If $d<mp$, then there exist non-negative integers $x$ and $y$ with
$x<n\le N<y$ such that $y-x\le mp$ and $a_x=a_y=1$.\\
\indent We now prove the claim. 
Note that $a_n=-1$ implies that $n \not\in S$ and $n-1 \in S$ and further we have $N \not \in S$ and $N-1 \in S$. Using $n-1\in S$ 
we have $n-1+mp\in S$. Since $n-1+mp\geq N$, $n-1+mp\in S$ and $N \not \in S$, it follows that $n-1+mp> N$.
So there is at least one $N<y<n+mp$ with $y-1\not\in S$ and $y\in S$, so $a_y=1$. But again $y-1\not \in S$ implies 
$y-1-mp\not \in S$ and we have $n-1\in S$. 
Therefore there exists an $x$ with $y-mp\leq x< n$ and $x-1\not\in S$ and $x\in S$, so $a_x=1$. Furthermore we 
have $y-x\leq mp$ and $x<n\leq N<y$.\\
{\tt Claim 2}: Under the conditions of Claim 1 we have $e_{j,x}+e_{j,y}\geq \min\{q-mp^2+p,p\}$.\\%p-(p^2-g)
\indent The proof is rather short. If $j-yp+1\geq p$, then $e_{j,y}=p$, 
since $$p+q-1-j+yp\geq (p+q-1-j+Np)+p \geq p,$$ and we are done. Otherwise $e_{j,y} = j-yp+1$ and
$$e_{j,x}\geq p+q-1-j+xp\geq p+q-1-j+yp-mp^2=p+q-mp^2-e_{j,y},$$ and hence $e_{j,x}+e_{j,y}\geq p+q-mp^2$.\\

In order to finish the proof of this lemma, we have to deal with the following two cases and show that
$d_j\ge -g$, respectively $d_j\ge -(p^2-g)$.\\
{\tt Case 1}: $0<g<p$.\\
If $d\ge fp$, then by (\ref{-g}) we have $d_j\ge -g$, so we may assume that $d<fp$. Using Claim 1 we 
find $x$ and $y$ as described in Claim 1. Now using Claim 2, we find that
\begin{eqnarray*}
d_j & = & \sum{a_k e_{j,k}}\geq e_{j,x}+\sum_{n\leq k \leq N}{a_k e_{j,k}}+e_{j,y}\\
& \geq & \min\{q-fp^2+p,p\}-p+p-\ldots+p-p = \min\{g+p,p\}-p = 0. 
\end{eqnarray*}
{\tt Case 2}: $p^2-p<g<p^2$.\\
If $d\geq (f+1)p$, then by (\ref{-g}) we have $d_j\ge p^2-g>0$, so we may assume that $d<(f+1)p$.
Now we can use again the Claims to find $x$ and $y$ as needed and we conclude that
\begin{eqnarray*}
d_j & = & \sum{a_k e_{j,k}}\geq e_{j,x}+\sum_{n\leq k \leq N}{a_k e_{j,k}}+e_{j,y}\\
& \geq & \min\{q-(f+1)p^2+p,p\}-p+p-\ldots+p-p\\
& =  & \min\{g+p-p^2,p\}-p = -(p^2-g),
\end{eqnarray*}
which finishes the proof.\qed\\

In the proof of the next lemma we use the notation ${\cal C}_{\le 0}(f)$ for 
${\cal C}(f)\cap {\mathbb Z}_{\le 0}$.
\begin{Lem}
\label{38}
Let $p$ and $q$ be distinct primes. Put $f_{38}=\Phi_p\Phi_{q}\Phi_{p^2q}$ and $\beta(p,q)=\min(p,q,q ({\rm mod~}p^2),p^2-q ({\rm mod~}p^2))$.
We have $${\cal C}(f_{38})=\cases{\{-2,0,1,2\} & if $q=2$;\cr
\{-1,1,2\} & if $p=2$ and $q=3$;\cr
[-\beta(p,q),\min(p,q)] & otherwise.}$$
\end{Lem}
{\it Proof}. Put $z_1=\min(p,q)$.
On noting that 
$$f_{38}\equiv \Phi_p\Phi_q \equiv {1\over (1-x)^2}\equiv \sum_{j=1}^{z_1}jx^{j-1}({\rm mod~}x^{z_1}),$$
we have $[1,\min(p,q)]\subseteq {\cal C}(f_{38})$.
This in combination with Lemma \ref{38flauw} shows that ${\cal C}_{>0}(f_{38})=[1,\min(p,q)]$. It remains 
to show that ${\cal C}_{\le 0}(f_{38})$ is as asserted in the statement of the lemma.\\
1) The case $q=2$. Here we have $\beta(p,2)=2$.\\
We have $\Phi_{2p^2}(x) = \Phi_{2p}(x^p) = \Phi_{p}(-x^p)$.
Then $f_{38}(x) = (1 + x + \ldots + x^{p-1}) (1 + x) (1 - x^p + x^{2p} - \ldots + x^{p^2-p}).$
Since
$(1 + x + \ldots + x^{p-1})(1 - x^p + x^{2p} - \ldots + x^{p^2-p}) =
1 + x + \ldots + x^{p-1} - x^p - x^{p+1} - \ldots - x^{2p-1} + \ldots + x^{p^2-p} + x^{p^2-p+1} 
+ \ldots + x^{p^2-1}$,
we have $f_{38}(x) = 1 + 2x + 2x^2 + \ldots + 2x^{p-1} - 2x^{p+1} - 2x^{p+2} - \cdots - 2x^{2p-1}
+ 2x^{2p+1} + \ldots + 2x^{p^2-1} + x^{p^2}$ and hence ${\cal C}(f_{38})=\{-2,0,1,2\}$.\\
2) The case $p=2$. Here we have $\beta(2,q)=1$.\\
If $q=3$ we have to show (cf. statement of this lemma) that 
${\cal C}_{\le 0}(f_{38})=\{-1\}$ (which follows by direct
calculation) and for $q\ge 5$ that ${\cal C}_{\le 0}(f_{38})=[-1,0]$.\\
\indent We have
\begin{eqnarray*}
	f_{38}(x) & = & \Phi_2(x)\Phi_q(x)\Phi_{2q}(x^2)\\
	& = & (1+2x+\ldots+2x^{q-1}+x^q)(1-x^2+x^4-x^6+\ldots+x^{2q-2}).
\end{eqnarray*}
Assume $q\ge 5$. It is easy to see that $d_3=0$. Furthermore, 
$$d_{q+1}=(-1)^{(q-1)/2}(1+\sum_{j=1}^{(q-1)/2}(-1)^j2)=-1.$$
It follows that ${\cal C}_{\le 0}(f_{38})=[-1,0]$.\\
3) The case where both $p$ and $q$ are odd.\\
Here we invoke Lemma \ref{38-1} and Lemma \ref{38-0}. \qed

\begin{Lem}
\label{39}
Let $p$ and $q$ be distinct primes. Put $f_{39}=(x-1)\Phi_p\Phi_{q}\Phi_{p^2q}$.
We have ${\cal C}(f_{39})=[-2,2]$.
\end{Lem}
{\it Proof}. Since $f_{39}=f_{36}(x^p-1)$ and $H(f_{36})=1$ by Lemma \ref{f36}, we immediately conclude that ${\cal C}(f_{39}) \subseteq [-2,2]$. Since $f_{39}$ is anti-self-reciprocal it is enough to show that, e.g. $0,1,2$ are in ${\cal C}(f_{39})$. Because $f_{39}$ is monic, this is clear for $1$. Write $f_{39}=\sum_j c_jx^j$. Note
that $f_{39}=\Psi_{pq}(x)\Phi_{pq}(x^p)$.\\
If $p>q$, we deduce from this and (\ref{invi}) that 
$$f_{39}(x)=(-1-x-\ldots-x^{q-1}+x^p+\ldots+x^{p+q-1})\Phi_{pq}(x^p).$$ 
We obtain 
$f_{39}(x)\equiv -1-x-\ldots-x^{q-1}+2x^p({\rm mod~}x^{p+1})$ on noting 
that $\Phi_{pq}(x^p)\equiv 1-x^p ({\rm mod~}x^{p+1})$. 
%It is also easy to note that there are only two combinations for $x^p$ which sum to a coefficient of two.
Hence $c_q=0$ and $c_p=2$.\\
If $p<q$, we find, using (\ref{invi}) again, 
\begin{equation}
\label{gutgut}
f_{39}(x)=(-1-x-\ldots-x^{p-1}+x^q+\ldots+x^{p+q-1})\Phi_{pq}(x^p)
\end{equation}
i) We start by showing that the coefficient 0 occurs in the $p<q$ case.\\

\noindent First assume $p=2$. \\
We have $\Phi_{pq}(x^p)=\Phi_{q}(-x^2)=1-x^2+x^4-\ldots+x^{2q-2}$. So if $q \equiv 1 ({\rm mod~}4)$, $[\Phi_{q}(-x^2)]_{x^{q-1}}=1$ and hence $c_q=1-1=0$. With $q \equiv 3 ({\rm mod~}4)$ we conclude $[\Phi_{q}(-x^2)]_{x^{q+1}}=1$ and hence $c_{q+1}=1-1=0$. So $c_qc_{q+1}=0$.\\
Next assume $p>2$. \\
If $q>2p$, then it is easy to see that $c_{2p}=0$.\\
In case $p<q<2p$ and $p \neq 3$, we have $c_{p^2-1}=0$.\\
The last case is $p=3$ and $q=5$. Here $0 \in {\cal C}(f_{39})$ follows by explicit 
computation ($c_{15}=c_{16}=0$).\\

\noindent ii) Next we show that the coefficients $-2$ (or $2$) also occur in the $p<q$ 
case (thus also $2$ (or $-2$) by self-reciprocity).\\

\noindent We write $\Phi_{pq}(x) = \sum{a_k}x^k$ and let $\rho$ and $\sigma$ be as in Lemma \ref{binary}.\\
We have $q=mp+g$ with $0<g<p$ and put $M=m+1$. 
We can write $M+k_Mpq=\rho_Mp+\sigma_Mq$ with $0\leq \rho_M<q$, $0\leq \sigma_M<p$ and $0\leq k_M\leq 1$. Note that $m<M<q$.\\
Now we study six different cases.\\
1) If $\rho_M\leq \rho$ and $\sigma_M \leq \sigma$ then $a_M=1$ by Lemma \ref{binary} and of course $a_1=-1$.
Now we determine the coefficient of $x^{p+q}$. Since $p+q=Mp+g$, we have  $[(-1-x-\ldots-x^{p-1})\Phi_{pq}(x^p)]_{x^{p+q}}=-a_M$.
Therefore $c_{p+q}=a_1-a_M=-2$.\\
\indent Before discussing the remaining five cases, we will establish the following observation.\\
{\tt Observation}: If $a_j=1$ and $a_{j+M}=-1$, then $c_{jp+p+q-1}=2$.\\
This is easy to check, we have $jp+p+q-1=(j+M)p+(g-1)$ and \mbox{$(jp+p+q-1)-q=jp+(p-1)$} with $0\leq g-1<p-1$. So 
$$[(-1-x-\ldots-x^{p-1})\Phi_{pq}(x^p)]_{x^{jp+p+q-1}}=-a_{j+M},$$ 
$$[(x^q+\ldots+x^{p+q-1})\Phi_{pq}(x^p)]_{x^{jp+p+q-1}}=[(1+\ldots+x^{p-1})\Phi_{pq}(x^p)]_{x^{jp+p-1}}=a_j,$$
resulting by (\ref{gutgut}) in $c_{jp+p+q-1}=-a_{j+M}+a_j=2$.\\
2) The second case we study is $\rho_M> \rho$ and $\sigma_M > \sigma$. Then $a_M=-1$ by 
Lemma \ref{binary} and since $a_0=1$, we can use the 
Observation to conclude that $c_{p+q-1}=2$.\\
3) The third case is $\rho_M> \rho$ and $\sigma_M =0$. But this does not arise, since otherwise it would 
follow (note that $k_M=0$ and $M<q$) that 
\begin{eqnarray*}
pq+1 & = & (\rho+1)p+(\sigma+1)q \leq (\rho_M-1+1)p+(\sigma+1)q\\
& = & M+(\sigma+1)q<(\sigma+2)q\leq pq,
\end{eqnarray*}
which is impossible.\\
4) The fourth case is $\rho_M> \rho$ and $0<\sigma_M \leq \sigma$.\\
Since $M<q$ and $\sigma_M>0$ we must have $k_M=1$. 
Put $j=(\sigma+1-\sigma_M)q$ and $j_1=(\sigma+1)q+\rho_Mp-pq$.
Since (use Lemma \ref{binary}) $$a_j=1,~a_{j_1}=-1,~j_1-j=\rho_Mp+\sigma_Mq-pq=M,$$ we can use the 
Observation to conclude that $c_{jp+p+q-1}=2$.\\
5) The case $\rho_M=0$ and $\sigma_M >\sigma$ does not arise, since otherwise we would
have $M+k_Mpq=\sigma_Mq<pq$ and 
therefore $k_M=0$ and so $M=\sigma_Mq\ge q$, contradicting $M<q$.\\
6) The last remaining case is $0<\rho_M\leq \rho$ and $\sigma_M > \sigma$.\\
Again, $M<q$ and $\sigma_M>\sigma>0$, so $k_M=1$. 
Put $j=(\rho+1-\rho_M)p$ and $j_1=(\rho+1)p+\sigma_Mq-pq$. 
Since (use Lemma \ref{binary}) $$a_j=1,~a_{j_1}=-1,~j_1-j=\rho_Mp+\sigma_Mq-pq=M,$$ we can use the Observation once again to conclude that $c_{jp+p+q-1}=2$.\qed

\subsubsection{The polynomials $f_{42}$ and $f_{43}$}
Let $p$ and $q$ be distinct primes. Put $f_{42}=\Phi_p\Phi_{pq}\Phi_{p^2q}=\sum c_jx^j$ and 
$f_{43}=(x-1)f_{42}=\sum d_jx^j$.
It is not difficult to find cases where only very few of the coefficients
of $f_{43}$ are equal to $2$. For example, if $(p,q)$ is in the following set:
 $$\{(11,241),(13,377),(17,577),(19,181),(29,421),
(41,3361),(43,3697)\},$$ there are precisely two coefficients equal to $2$ (as computed by
Yves Gallot). This suggests that perhaps the following results is not so easy to establish.
\begin{Lem}
\label{T43}
We have $${\cal C}(f_{43})=\cases{\{-2,-1,1,2\} & if $q=2$;\cr 
[-2,2] & otherwise.}$$
\end{Lem}
The analogue of this result for $f_{42}$ is easy enough. Note that 
$${\rm deg}(f_{42})=p^2(q-1)+p-q.$$
\begin{Lem}
\label{42}
We have
${\cal C}(f_{42})=[-1,1]$.
\end{Lem}
{\it Proof}. Write $f_{42}=\sum_j c_jx^j$. Note that
$$f_{42}={(x^p-1)(x^{p^2q}-1)\over (x^q-1)(x^{p^2}-1)}.$$
Around $x=0$, $f_{42}$ has power series
\begin{equation}
\label{taylor}
(1+x^q+x^{2q}+\cdots)(1-x^p+x^{p^2}-x^{p^2+p}+\cdots+x^{(q-1)p^2}-x^{(q-1)p^2+p}).
\end{equation}
Note that if $c_j\ge 2$, then there exist non-negative
$\alpha_1,\alpha_2,\beta_1$ and $\beta_2$ such that
$$\alpha_1\ne \alpha_2,~\beta_1\ne \beta_2,~j=\alpha_1q+\beta_1p^2=\alpha_2q+\beta_2p^2
\le {\rm deg}(f_{42})<p^2q.$$
This is impossible. By a similar argument one sees that $c_j\ge -1$. Since clearly
$[-1,1]\subseteq {\cal C}(f_{42})$, the proof is completed.\qed\\

\noindent Indeed, some work needs to be done to infer that $\{-2,2\}\subseteq {\cal C}(f_{43})$. The idea is to
show that in $f_{42}$ the combinations $1,-1$ and $-1,1$
appear as consecutive coefficients and then use that $f_{43}=(x-1)f_{42}$.\\
\indent Let us denote by $\{a;b\}$ the smallest non-negative integer $m$ such that $m\equiv a({\rm mod~}b)$.
\begin{Lem}
\label{43}
Write $f_{42}=\sum c_jx^j$ and $f_{43}=\sum d_jx^j$.
Put 
$$k_1=1+\{{p-1\over p^2};q\}p^2{\rm ~and~}k_2=1+\{{p-1\over q};p^2\}q.$$
{\rm 1)} Suppose that $1<k_1\le {\rm deg}(f_{42})$.
If furthermore, 
\begin{equation}
\label{r1}
\{{1\over q};p^2\}q+\{{1\over p};q\}p^2>p^2q
\end{equation}
and 
\begin{equation}
\label{s1}
\{{-1\over q};p\}pq+\{{-1\over p^2};q\}p^2+p+1>p^2q,
\end{equation}
then $c_{k_1-1}=1$, $c_{k_1}=-1$ and $d_{k_1}=2$.\\
{\rm 2)} Suppose that $k_2\le {\rm deg}(f_{42})$.
If furthermore,
\begin{equation}
\label{r2}
\{{-1\over q};p^2\}q+\{{-1\over p};q\}p^2+p+1>p^2q
\end{equation}
and 
\begin{equation}
\label{s2}
\{{1\over q};p\}pq+\{{1\over p^2};q\}p^2>p^2q,
\end{equation}
then $c_{k_2-1}=1$, $c_{k_2}=-1$ and $d_{k_2}=2$.
\end{Lem}
{\it Proof}. We say that $k$ is {\it $p$-representable} if we can write 
$k=m_1q+m_2p^2$ with $m_1\ge 0$ and $0\le m_2\le q-1$. We say that 
$k$ is {\it $m$-representable} if we can write $k=n_1q+n_2p^2+p$ with $n_1\ge 0$ and
$0\le n_2\le q-1$. {}From the proof of Lemma \ref{42} it follows that if $k\le {\rm deg}(f_{42})$,
then $k$ can be $p$-representable in at most one way and be $m$-representable in at most
one way. {}From this and (\ref{taylor}), we infer that if $k\le {\rm deg}(f_{42})$, then 
\begin{equation}
\label{minnie}
c_k=\cases{1 & if $k$ is $p$-representable, but not $m$-representable;\cr
-1 & if $k$ is $m$-representable, but not $p$-representable;\cr
0 & otherwise.}
\end{equation}
We have
\begin{equation}
\label{dubbel}
\cases{k_1\equiv 1({\rm mod~}p^2);\cr k_1\equiv p({\rm mod~}q),} {\rm ~and~}
\cases{k_2\equiv p({\rm mod~}p^2);\cr k_2\equiv 1({\rm mod~}q).}
\end{equation}
Suppose that $k_1\le {\rm deg}(f_{42})$.
Clearly $k_1$ is $m$-representable, because $k_1>1$ implies $k_1>p$. Condition (\ref{r1}) ensures that $k_1$ is not
$p$-representable. Thus, by (\ref{taylor}), we have $c_{k_1}=-1$. On the other hand
we see that $k_1-1$ is $p$-representable, but not $m$-representable by (\ref{s1}). It
follows that $c_{k_1-1}=1$. Since $d_{k_1}=c_{k_1-1}-c_{k_1}=1-(-1)=2$, we have established
part 1.  Part 2 can be derived in a similar way, but here it is not needed to require $k_2>1$. \qed\\

\noindent We will show that some of the numbers appearing in the latter lemma are actually equal. For this
the reciprocity law formulated in Corollary \ref{recipro} is needed. 
As usual by $(m,n)$ we denote the greatest common divisor of $m$ and $n$.
\begin{Lem} 
\label{andreas}
Suppose that $a>1$ and $b>1$ are coprime integers. 
Then $$(a-\{\frac 1b;a\},\{\frac 1a;b\})=(\{\frac 1b;a\},b-\{\frac 1a;b\})=1.$$
\end{Lem}
\begin{cor}
\label{recipro}
Suppose that both $a>1$ and $b>1$ are odd and coprime. Then the congruence $\{{1\over a};b\}\equiv \{{1\over b};a\}
({\rm mod~}2)$ holds.
\end{cor}
{\it Proof}. If $\{{1\over a};b\}$ is even,  then $a-\{\frac 1b;a\}$ must be odd and hence 
$\{\frac 1b;a\}$ is even. If $\{{1\over a};b\}$ is odd,  then $b-\{\frac 1a;b\}$ is even and hence 
$\{\frac 1b;a\}$ must be odd. \qed\\

\noindent {\it Proof of Lemma} \ref{andreas}. Put 
$\delta(a,b)=(\{\frac 1a;b\})(\{\frac 1b;a\})-(a-\{\frac 1b;a\})(b-\{\frac 1a;b\}).$
It is enough to show that $\delta(a,b)=1$.
Since clearly
$-ab+1<\delta(a,b)<ab,$ it is enough to show that 
$\delta(a,b) \equiv 1 ({\rm mod~}ab)$.
We have 
$\delta(a,b) \equiv \{\frac 1b;a\}b \equiv 1({\rm mod~}a)$
and $\delta(a,b) \equiv 
\{\frac 1a;b\}a \equiv 1 ({\rm mod~}b),$ and on invoking the Chinese remainder theorem the proof is
completed.\qed

\begin{Lem}
\label{decker}
We have
$$\{\frac 1q;p^2\}q+\{\frac 1p;q\}p^2=\{\frac{-1}q;p\}pq+\{\frac{-1}{p^2};q\}p^2+p+1$$
and
$$\{\frac{-1}{q};p^2\}q+\{\frac{-1}{p};q\}p^2+p+1=\{\frac 1q;p\}pq+\{\frac 1{p^2};q\}p^2.$$
\end{Lem}
{\it Proof}. Denote the numbers appearing in the left hand sides of {\rm (\ref{r1}), (\ref{s1}), (\ref{r2})}
and {\rm (\ref{s2})}, by $r_1(p,q),s_1(p,q),r_2(p,q),s_2(p,q)$, respectively. 
We have to show that $r_1(p,q)=s_1(p,q)$ and $r_2(p,q)=s_2(p,q)$. On noting that
$\{{-1\over q};p\}=p-\{{1\over q};p\}$, etc., it is easily seen that $r_1(p,q)=s_1(p,q)$  implies
$r_2(p,q)=s_2(p,q)$, thus it is enough to show that $r_1(p,q)=s_1(p,q)$. 
By considering  $r_1,r_2,s_1,s_2$ modulo $p^2$ and
$q$ and invoking the Chinese remainder theorem we infer that
\begin{equation}
\label{boeha}
k_j\equiv r_j(p,q)\equiv s_j(p,q)({\rm mod~}p^2q)~{\rm for~}1\le j\le 2.
\end{equation}
Note that 
\begin{equation}
\label{ers}
\{r_j(p,q),s_j(p,q)\}\subseteq \{k_j,k_j+p^2q\} ~{\rm ~for~}1\le j\le 2.
\end{equation} 
Thus it suffices to establish that $r_1(p,q)\equiv s_1(p,q)({\rm mod~}2p^2q)$ in order to show 
that $r_1(p,q)=s_1(p,q)$.\\
1) $p=2$. Recall that the Legendre symbol $({-1\over q})$ equals $(-1)^{(q-1)/2}$ in case $q$ is odd.
We have $r_1(2,q)=\{\frac 1q;4\}q+\{\frac 12;q\}4=4q+2-({-1\over q})q$, on noting that
$\{\frac 1q;4\}=2-({-1\over q})$ and $\{\frac 12;q\}=(q+1)/2$. On noting that
$\{\frac{-1}q;2\}=1$ and $\{\frac{-1}{4};q\}4=(2-({-1\over q}))q-1$, one infers that
$$s_1(2,q)=\{\frac{-1}q;2\}2q+\{\frac{-1}{4};q\}4+2+1=4q+2-({-1\over q})q=r_1(2,q).$$
2) $q=2$. By an argument easier than that for case 1 one infers that $r_1(p,2)=s_1(p,2)=2p^2+1$.\\
3) $p,q$ odd. It suffices to show that $r_1(p,q)\equiv s_1(p,q)({\rm mod~}2)$. Now
using Corollary \ref{recipro} we have $\{\frac 1q;p^2\}=\{\frac 1{p^2};q\}({\rm mod~}2)$
and $\{\frac 1p;q\}=\{\frac 1q;p\}({\rm mod~}2)$ and hence
\begin{eqnarray}
\{\frac 1q;p^2\}q+\{\frac 1p;q\}p^2 & \equiv & \{\frac 1q;p^2\}+\{\frac 1p;q\} \equiv \{\frac{1}{p^2};q\}+\{\frac 1q;p\}\nonumber\\
& \equiv & q-\{\frac{1}{p^2};q\}+p-\{\frac 1q;p\}  \equiv  \{\frac{-1}{p^2};q\}+\{\frac{-1}q;p\}\nonumber\\
& \equiv & \{\frac{-1}q;p\}pq+\{\frac{-1}{p^2};q\}p^2+p+1~({\rm mod~}2),\nonumber
\end{eqnarray}
which finishes the proof.\qed 

\begin{Lem}
\label{43a}
Write $f_{43}=\sum_j d_jx^j$.
There is a unique integer $1\le j\le 2$ such that the conditions of part $j$ of 
Lemma {\rm \ref{43}} are satisfied and hence $d_{k_j}=2$. 
Furthermore, $d_{{\rm deg}(f_{42})-k_j+1}=-2$.
\end{Lem}
{\it Proof}. We consider the cases  $p\not\equiv 1({\rm mod~}q)$ and $p\equiv 1({\rm mod~}q)$ separately.\\
i) The case $p\not\equiv 1({\rm mod~}q)$.\\ 
We have $q\ge 3$ and  $k_1>1$.
{}From (\ref{dubbel}) we infer that $k_1+k_2\equiv 1+p({\rm mod~}p^2q)$. Since 
clearly $1+p<1+p^2\le k_1+k_2<1+p+2p^2q$,
we infer that
\begin{equation}
\label{k1k2}
k_1+k_2=1+p+p^2q.
\end{equation} 
Let us suppose that $k_1\ge p^2(q-1)+p-q+1={\rm deg}(f_{43})=1+{\rm deg}(f_{42})$. By (\ref{k1k2}) we then have
$k_2\le p^2+q$. Since $q \geq 3$ and $p^2+p\geq 6$ it follows that
\begin{eqnarray*}
k_2 & \leq & p^2+q\leq 2p^2+p+q-6 = 3(p^2-2)+p+q-p^2\\
& \leq & q(p^2-2)+p+q-p^2 = qp^2+p-q-p^2 = (q-1)p^2+p-q,
\end{eqnarray*}
so $k_2\le {\rm \deg}(f_{42})$. Since $r_2(p,q)>p^2+q\ge k_2$ and
$r_2(p,q)\equiv k_2({\rm mod~}p^2q)$, we have 
$r_2(p,q)=k_2+p^2q>p^2q$. Since $r_2(p,q)=s_2(p,q)$ by Lemma \ref{decker}, it follows
that if $k_1>{\rm \deg}(f_{42})$ and thus the conditions of part 1 (of Lemma \ref{43}) are not satisfied, then the
conditions of part 2 are satisfied. By a similar argument we infer that if 
$k_2>{\rm \deg}(f_{42})$ and thus the conditions of part 2 are not satisfied, then the
conditions of part 1 are satisfied.\\
\indent It remains to deal with the case where $k_j\le {\rm \deg}(f_{42})$ for
$1\le j\le 2$. Note that
\begin{equation}
\label{reenrtwee}
r_1(p,q)+r_2(p,q)=1+p+2p^2q.
\end{equation}
Hence $r_j(p,q)>p^2q$ for some $1\le j\le 2$. Let us assume w.l.o.g. that
$r_2(p,q)>p^2q$. Now if $r_1(p,q)>p^2q$, then on using (\ref{ers}) we find
$$r_1(p,q)+r_2(p,q)=k_1+p^2q+r_2(p,q)>1+p^2+2p^2q,$$
contradicting (\ref{reenrtwee}) and hence the conditions of both part 1 and part 2 cannot be
satisfied at the same time.\\
ii) The case $p\equiv 1({\rm mod~}q)$.\\ 
Here we can write $p=kq+1$, with $k\ge 1$. We have $k_1=1+0p^2=1$ and 
$k_2=1+kq=p$ and hence the conditions of part 1 are not satisfied. 
We have to show that the conditions of part 2 are satisfied.
Obviously $k_2=p\leq {\rm \deg}(f_{42})$. On noting that $\{{-1\over q};p^2\}=k^2q+2k$ and $\{{-1\over p};q\}=q-1$, 
the left side of equation (\ref{r2}) becomes:
$$(k^2q+2k)q+(q-1)p^2+p+1=k^2q^2+2kq+p^2q-p^2+p+1=p^2q+p>p^2q.$$
Similarly we have for the left side of equation (\ref{s2}):
$$(p-k)pq+1\cdot p^2=p^2q-p(p-1)+p^2=p^2q+p>p^2q.$$
(Alternatively one can invoke Lemma \ref{decker} to deduce that the left hand side of 
(\ref{s2}) equals the left hand side of (\ref{r2}) and hence exceeds $p^2q$.)\\
\indent In both cases i) and ii), we conclude that there is a unique integer $j$ such the conditions
of part j of Lemma \ref{43} are satisfied.\\
\indent The final assertion follows on noting that $f_{42}$ is self-reciprocal and using that 
$f_{43}=(x-1)f_{42}$.\qed\\
   
\noindent {\tt Examples}. Using the latter lemma, one can derive the following examples.\\
1) If $p^2+p-1\equiv 0({\rm mod~}q)$, then $d_{p^2+p}=2$.\\
2) If $p^2-p+1\equiv 0({\rm mod~}q)$, then $d_{p^2+1}=2$.\\
3) If $q\equiv 1-p({\rm mod~}p^2)$, then $d_{q+p}=2$.\\
4) If $q\equiv p-1({\rm mod~}p^2)$, then $d_{q+1}=2$.\\
5) If $p\equiv 1({\rm mod~}q)$, then $d_p=2$.\\

\noindent {\it Proof of Lemma} \ref{T43}. By (\ref{zon}) and Lemma \ref{42} we find that
${\cal C}(f_{43})\subseteq [-2,2]$.  By Lemma \ref{43a} we have $\{-2,2\}\subseteq C(f_{43})$.
Since $d_0=-1$ and $d_{{\rm deg}(f_{43})}=1$, it remains to be determined when $0\in {\cal C}(f_{43})$.
If both $p$ and $q$ are odd, then $d_2=0$. If $q=2$, then  $f_{43}$ has the power series (around $x=0$)
$$f_{43}=(-1+x^p-x^{p^2}+x^{p^2+p})(1-x+x^2-x^3+x^4-x^5+\cdots )$$
and since $p$ is odd we find that
$d_j\ne 0$ for $j\le {\rm deg}(f_{43})=p^2+p-1$ and hence 
$0\not\in {\cal C}(f_{43})$.\\ 
If $p=2$, then $f_{43}$ has the power series (around $x=0$)
$$f_{43}=(1+x^q+x^{2q}+x^{3q})\sum_{k=0}^{\infty}(-x^{4k}+x^{4k+1}+x^{4k+2}-x^{4k+3}).$$
{}From this we see that $d_q=0$ if $q\equiv 1({\rm mod~}4)$ and $d_{q+1}=0$ if
$q\equiv 3({\rm mod~}4)$. Since $q+1<{\rm deg}(f_{43})=3q-1$, it follows that $0\in {\cal C}(f_{43})$ if $p=2$.\qed

\section{The proof of the main theorem}
\label{proef}
{\it Proof of Theorem} \ref{MAIN}. {}From $x^n-1=\prod_{d|n}\Phi_d(x)$ and the fact that the
$\Phi_d$ are irreducible over the rationals, we infer that any divisor of $x^n-1$ with integer
coefficients is of the form $\pm \prod_{d|n}\Phi_d^{e_d}(x)$, with $e_i\in \{0,1\}$.
Thus we have $2^{d(n)}$ monic divisors, where $d(n)$ denotes the number of divisors of $n$.\\
\indent {}From the identity
\begin{equation}
\label{blu}
x^{p^2q}-1=\Phi_1(x)\Phi_p(x)\Phi_q(x)\Phi_{pq}(x)\Phi_{p^2}(x)\Phi_{p^2q}(x),
\end{equation}
we infer that $x^{p^2q}-1$ has 64 divisors. We denote these by $f_0,\ldots,f_{63}$.
If $k=\sum_{j=0}^5 k_j2^j$ is the base 2 expansion of $k$, then we put
$$f_k(x)=\Phi_1(x)^{k_0}\Phi_p(x)^{k_1}\Phi_q(x)^{k_2}\Phi_{pq}(x)^{k_3}\Phi_{p^2}(x)^{k_4}
\Phi_{p^2q}(x)^{k_5}.$$ Thus $\{f_0(x),\ldots,f_{63}(x)\}$ is the set of all 
monic divisors of $x^{p^2q}-1$.
Note that $\Phi_1(x)=x-1$,
$\Phi_p(x)=1+x+\cdots+x^{p-1}$ and $\Phi_q(x)=1+x+\cdots+x^{q-1}$. Thus these three
divisors have all height 1. By Lemma \ref{binary}
we have $H(\Phi_{pq}(x))=1$.
On noting that $\Phi_{p^2}(x)=\Phi_p(x^p)$ and $\Phi_{p^2q}(x)=\Phi_{pq}(x^p)$, it then
follows that each of the six cyclotomic polynomials appearing in (\ref{blu}) is flat.\\
\indent We will only establish the less trivial cases in Table 1, the easier ones being left
as exercises to the reader. (Note that for some polynomials like $f_{19}$ we have given more than
one argument.)\\
-$f_0,f_1,f_2,f_3,f_4,f_5,f_{16},f_{17},f_{18},f_{19}$: Use Theorem \ref{riant}.\\
-$f_6$. Use Lemma \ref{vijf}.\\
-$f_7$. Use Corollary \ref{ojoj}.\\
-$f_8$: Use Lemma \ref{tau}.\\
-$f_9$. Use Lemma \ref{nulgat}.\\
-$f_{16},f_{17},f_{18},f_{32},f_{33},f_{34}$: Use identity (\ref{bier}).\\
-$f_{20}$: See Lemma \ref{20}.\\
-$f_{21},f_{37}$: Note that $\Phi_1(x)\Phi_q(x)=x^q-1$.\\ 
-$f_{22}$: See Corollary \ref{22}.\\
-$f_{19}$, $f_{23}, f_{27}$. Use that $\Phi_1(x)\Phi_p(x)\Phi_{p^2}(x)=x^{p^2}-1$.\\
-$f_{24}$: Invoke Lemma \ref{koffie}.\\
-$f_{25}$: Invoke Lemma \ref{25}.\\
-$f_{26},f_{27}$: Invoke Lemma \ref{26}.\\
-$f_{28}$: We have $f_{28}=\Phi_p(x^p)\Phi_q(x^p)$. On invoking the result that
${\cal C}(\Phi_p\Phi_q)=[1,\min(p,q)]$ (follows by Lemma \ref{vijf}), the assertion follows.\\
-$f_{29}$: If $p=2$, then consecutive coefficients in $f_{28}$ are distinct and hence $0\not\in {\cal
C}(f_{29})$.\\
-$f_{30}$: See Lemma \ref{30}.\\
-$f_{31}$. Note that $f_{31}=(x^{p^2q}-1)/\Phi_{p^2q}(x)=\Psi_{p^2q}(x)=\Psi_{pq}(x^p)$. Thus,
${\cal C}(f_{31})=[-1,1]$ by (\ref{invi}).\\
-$f_{34}$: Using (\ref{flauuw}) we find that ${\cal C}(f_{34})={\cal C}(f_8)$.\\
-$f_{35}$: $f_{35}=(x^p-1)\Phi_{pq}(x^p)=f_9(x^p)$. 
It follows that ${\cal C}(f_{35})={\cal C}(f_9)\cup \{0\}$.
Now invoke Lemma \ref{nulgat}.\\
-$f_{36}$: Invoke Lemma \ref{f36}.\\
-$f_{37}$: We have $f_{37}=(x^q-1)\Phi_{pq}(x^p)$. Noting that $q+jp\neq kp$, we infer 
that ${\cal C}(f_{37})=[-1,1]$.\\
-$f_{38}$: Invoke Lemma \ref{38}.\\
-$f_{39}$. Invoke Lemma \ref{39}.\\
-$f_{40}$: We have $f_{40}=\tau_{p^2,q}(x)$. Now invoke Lemma \ref{tau}.\\
-$f_{41}$. Invoke Lemma \ref{nulgat}.\\
-$f_{42}$. Invoke Lemma \ref{42}.\\
-$f_{43}$. Invoke Lemma \ref{T43}.\\
-$f_{44}$: We have $\Phi_q(x)\Phi_{pq}(x)\Phi_{pq}(x^p)=(x^{p^2q}-1)/(x^{p^2}-1)$.\\
\medskip
-$f_{48},f_{49},\ldots,f_{63}$.\\
Let $0\le j\le 15$.
Note that
$$f_{j+48}=f_j\Phi_{p^2}(x)\Phi_{p^2q}(x)=
f_j\Phi_{p}(x^p)\Phi_{pq}(x^p)=f_j(1+x^{pq}+x^{2pq}+\cdots+x^{(p-1)pq}),$$
 it follows by (\ref{bier}) that if deg$(f_j)<pq-1$, then ${\cal C}(f_{j+48})={\cal C}(f_j)\cup \{0\}$.\\
We have deg$(f_j)\geq pq-1$ iff\\
-$q=2$, $j=11$;\\
-$p=2$, $j=13$;\\
-$j=14$;\\
-$j=15$.\\
Using these two observations and Table 1A, one easily arrives at Table 1D. \qed

\medskip

%\vfil\eject
\centerline{{\bf Table 1}}
\medskip
\noindent Table 1 comes in 4 parts, 1A, 1B, 1C and 1D, each listing ${\cal C}(f)$ for 16
monic divisors of $x^{p^2q}-1$. For each of the tables there are some exceptions
to the set ${\cal C}(f)$ given in the table and these are listed directly below the 
table. If $\min(p,q)>3$, then there are no exceptions and ${\cal C}(f)$ can be read
of directly from the table.\\

\centerline{{\bf Table 1A}}
\medskip
\begin{center}
\begin{tabular}{|c|c|c|c|c|c|c|c|c|}
\hline
$f$ & $\Phi_1(x)$ & $\Phi_p(x)$ & $\Phi_q(x)$ & $\Phi_{pq}(x)$ & $\Phi_{p^2}(x)$ & $\Phi_{p^2q}(x)$ &
${\cal C}(f)$\\
\hline
\hline
0 & 0 & 0 & 0 & 0 & 0 & 0 & $\{1\}$\\ 
\hline
1 & 1 & 0 & 0 & 0 & 0 & 0 & $\{-1,1\}$\\ 
\hline
2 & 0 & 1 & 0 & 0 & 0 & 0 & $\{1\}$\\ 
\hline
3 & 1 & 1 & 0 & 0 & 0 & 0 & $[-1,1]$\\ 
\hline
4 & 0 & 0 & 1 & 0 & 0 & 0 & $\{1\}$\\ 
\hline
5 & 1 & 0 & 1 & 0 & 0 & 0 & $[-1,1]$\\ 
\hline
6 & 0 & 1 & 1 & 0 & 0 & 0 & $[1,\min(p,q)]$ \\ 
\hline
7 & 1 & 1 & 1 & 0 & 0 & 0 & $[-1,1]$\\ 
\hline
8 & 0 & 0 & 0 & 1 & 0 & 0 & $[-1,1]$\\ 
\hline
9 & 1 & 0 & 0 & 1 & 0 & 0 & $[-2,2]$\\ 
\hline
10 & 0 & 1 & 0 & 1 & 0 & 0 & $[0,1]$\\ 
\hline
11 & 1 & 1 & 0 & 1 & 0 & 0 & $[-1,1]$\\ 
\hline
12 & 0 & 0 & 1 & 1 & 0 & 0 & $[0,1]$\\ 
\hline
13 & 1 & 0 & 1 & 1 & 0 & 0 & $[-1,1]$\\ 
\hline
14 & 0 & 1 & 1 & 1 & 0 & 0 & $\{1\}$\\ 
\hline
15 & 1 & 1 & 1 & 1 & 0 & 0 & $[-1,1]$\\ 
\hline
\end{tabular}
\end{center}
\vfil\eject
\noindent If $\min(p,q)=2$, then ${\cal C}(f_8)=\{-1,1\}$.\\
If $\min(p,q)\le 3$, then ${\cal C}(f_9)=\{-2,-1,1,2\}$.\\
If $q=2$, then ${\cal C}(f_{11})=\{-1,1\}$.\\
If $p=2$, then ${\cal C}(f_{13})=\{-1,1\}$.\\

\noindent We put $\alpha(p,q)=\min([{q-1\over p}]+1,p)$.\\
By $p^*$ we denote the unique integer with $1\le p^*<q$ such that
$pp^*\equiv 1({\rm mod~}q)$.\\
We define $\gamma(p,q)=\min(p,p^*)+\min(p,q-p^*)$.\\

%\vfil\eject
\centerline{{\bf Table 1B}}
\begin{center}
\begin{tabular}{|c|c|c|c|c|c|c|c|c|}
\hline
$f$ & $\Phi_1(x)$ & $\Phi_p(x)$ & $\Phi_q(x)$ & $\Phi_{pq}(x)$ & $\Phi_{p^2}(x)$ & 
$\Phi_{p^2q}(x)$ & ${\cal C}(f)$\\
\hline
\hline
16 &0 & 0 & 0 & 0 & 1 & 0 & $[0,1]$\\ 
\hline
17 & 1 & 0 & 0 & 0 & 1 & 0 & $[-1,1]$\\ 
\hline
18 & 0 & 1 & 0 & 0 & 1 & 0 & $\{1\}$\\ 
\hline
19 & 1 & 1 & 0 & 0 & 1 & 0 & $[-1,1]$\\ 
\hline
20 & 0 & 0 & 1 & 0 & 1 & 0 & $[\min([{q\over p}],1),\alpha(p,q)]$\\ 
\hline
21 & 1 & 0 & 1 & 0 & 1 & 0 & $[-1,1]$\\ 
\hline
22 & 0 & 1 & 1 & 0 & 1 & 0 & $[1,\min(p^2,q)]$\\ 
\hline
23 & 1 & 1 & 1 & 0 & 1 & 0 & $[-1,1]$\\ 
\hline
24 & 0 & 0 & 0 & 1 & 1 & 0 & $[-\min(p,q-p^*),\min(p,p^*)]$\\ 
\hline
25 & 1 & 0 & 0 & 1 & 1 & 0 & $[-\gamma(p,q),\gamma(p,q)]$\\ 
\hline
26 & 0 & 1 & 0 & 1 & 1 & 0 & $[0,1]$\\ 
\hline
27 & 1 & 1 & 0 & 1 & 1 & 0 & $[-1,1]$\\ 
\hline
28 & 0 & 0 & 1 & 1 & 1 & 0 & $[0,\min(p,q)]$\\ 
\hline
29 & 1 & 0 & 1 & 1 & 1 & 0 & $[-\min(p,q),\min(p,q)]$\\ 
\hline
30 & 0 & 1 & 1 & 1 & 1 & 0 & $[1,\min(p,q)]$\\ 
\hline
31 & 1 & 1 & 1 & 1 & 1 & 0 & $[-1,1]$ \\ 
\hline
\end{tabular}
\end{center}
If $p=2$, then ${\cal C}(f_{17})=\{-1,1\}$.\\
If $\min(p,q)=2$, then ${\cal C}(f_{24})=[-\min(p,q-p^*),\min(p,p^*)]\backslash\{0\}$.\\
If $p\le 3$ and $q\neq 2$, then ${\cal C}(f_{25})=[-\gamma(p,q),\gamma(p,q)]\backslash\{0\}$.\\
If $p=2$, then ${\cal C}(f_{29})=\{-2,-1,1,2\}=[-\min(2,q),\min(2,q)]\backslash \{0\}$.\\
\vfil\eject
\centerline{{\bf Table 1C}}
\begin{center}
\begin{tabular}{|c|c|c|c|c|c|c|c|c|}
\hline
$f$ & $\Phi_1(x)$ & $\Phi_p(x)$ & $\Phi_q(x)$ & $\Phi_{pq}(x)$ & $\Phi_{p^2}(x)$ & 
$\Phi_{p^2q}(x)$ & ${\cal C}(f)$ \\
\hline
\hline
32 & 0 & 0 & 0 & 0 & 0 & 1 & $[-1,1]$\\ 
\hline
33 & 1 & 0 & 0 & 0 & 0 & 1 & $[-1,1]$\\ 
\hline
34 & 0 & 1 & 0 & 0 & 0 & 1 & $[-1,1]$\\ 
\hline
35 & 1 & 1 & 0 & 0 & 0 & 1 & $[-2,2]$\\ 
\hline
36 & 0 & 0 & 1 & 0 & 0 & 1 & $[-1,1]$ \\ 
\hline
37 & 1 & 0 & 1 & 0 & 0 & 1 & $[-1,1]$ \\ 
\hline
38 & 0 & 1 & 1 & 0 & 0 & 1 & $[-\beta(p,q),\min(p,q)]$\\ 
\hline
39 & 1 & 1 & 1 & 0 & 0 & 1 & $[-2,2]$\\ 
\hline
40 & 0 & 0 & 0 & 1 & 0 & 1 & $[-1,1]$\\ 
\hline
41 & 1 & 0 & 0 & 1 & 0 & 1 & $[-2,2]$\\ 
\hline
42 & 0 & 1 & 0 & 1 & 0 & 1 & $[-1,1]$\\ 
\hline
43 & 1 & 1 & 0 & 1 & 0 & 1 & $[-2,2]$\\ 
\hline
44 & 0 & 0 & 1 & 1 & 0 & 1 & $[0,1]$\\ 
\hline
45 & 1 & 0 & 1 & 1 & 0 & 1 & $[-1,1]$\\ 
\hline
46 & 0 & 1 & 1 & 1 & 0 & 1 & $[0,1]$\\ 
\hline
47 & 1 & 1 & 1 & 1 & 0 & 1 & $[-1,1]$\\ 
\hline
\end{tabular}
\end{center}
%\vfil\eject
\noindent We put $\beta(p,q)=\min(p,q,q ({\rm mod~}p^2),p^2-q ({\rm mod~}p^2))$.\\
If $p=2$, then ${\cal C}(f_{33})=\{-1,1\}$.\\
If $\min(p,q)=2$, then ${\cal C}(f_{34})=\{-1,1\}$.\\
If $q=2$, then ${\cal C}(f_{38})=\{-2,0,1,2\}$.\\
If $q=3$ and $p=2$, then ${\cal C}(f_{38})=\{-1,1,2\}$.\\
If $q=2$, then ${\cal C}(f_{40})=\{-1,1\}$.\\
If $q\le 3$, then ${\cal C}(f_{41})=\{-2,-1,1,2\}$.\\
If $q=2$, then ${\cal C}(f_{43})=\{-2,-1,1,2\}$.\\

%\vfil\eject
\centerline{{\bf Table 1D}}
\begin{center}
\begin{tabular}{|c|c|c|c|c|c|c|c|c|}
\hline
$f$ & $\Phi_1(x)$ & $\Phi_p(x)$ & $\Phi_q(x)$ & $\Phi_{pq}(x)$ & $\Phi_{p^2}(x)$ & $\Phi_{p^2q}(x)$ &
${\cal C}(f)$\\
\hline
\hline
48 & 0 & 0 & 0 & 0 & 1 & 1 & $[0,1]$\\ 
\hline
49 & 1 & 0 & 0 & 0 & 1 & 1 & $[-1,1]$\\ 
\hline
50 & 0 & 1 & 0 & 0 & 1 & 1 & $[0,1]$\\ 
\hline
51 & 1 & 1 & 0 & 0 & 1 & 1 & $[-1,1]$\\ 
\hline
52 & 0 & 0 & 1 & 0 & 1 & 1 & $[0,1]$\\ 
\hline
53 & 1 & 0 & 1 & 0 & 1 & 1 & $[-1,1]$\\ 
\hline
54 & 0 & 1 & 1 & 0 & 1 & 1 & $[0,\min(p,q)]$ \\ 
\hline
55 & 1 & 1 & 1 & 0 & 1 & 1 & $[-1,1]$\\ 
\hline
56 & 0 & 0 & 0 & 1 & 1 & 1 & $[-1,1]$\\ 
\hline
57 & 1 & 0 & 0 & 1 & 1 & 1 & $[-2,2]$\\ 
\hline
58 & 0 & 1 & 0 & 1 & 1 & 1 & $[0,1]$\\ 
\hline
59 & 1 & 1 & 0 & 1 & 1 & 1 & $[-1,1]$\\ 
\hline
60 & 0 & 0 & 1 & 1 & 1 & 1 & $[0,1]$\\ 
\hline
61 & 1 & 0 & 1 & 1 & 1 & 1 & $[-1,1]$\\ 
\hline
62 & 0 & 1 & 1 & 1 & 1 & 1 & $\{1\}$\\ 
\hline
63 & 1 & 1 & 1 & 1 & 1 & 1 & $[-1,1]$\\ 
\hline
\end{tabular}
\end{center}
If $q=2$, then ${\cal C}(f_{59})=\{-1,1\}$.\\
If $p=2$, then ${\cal C}(f_{61})=\{-1,1\}$.

\subsection{Compact reformulation of Theorem {\rm \ref{MAIN}}}
\label{compact}
For reference purposes a more compact version of Theorem \ref{MAIN} might be useful.
We give it here (this reformulation was given by Yves Gallot).
\begin{Thm}
Let $p$ and $q$ be distinct primes. Let $f(x) \in \mathbb{Z}[x]$
be a monic divisor of $x^{p^2 q} - 1$.
There exists an integer $k = \sum_{j=0}^{5} k_j 2^j$
with $k_j \in \{0, 1 \,\}$ (the binary expansion of $k$) such that
$$f(x) = f_k(x) = \Phi_1^{k_0}\cdot \Phi_p^{k_1}\cdot \Phi_q^{k_2}\cdot
  \Phi_{pq}^{k_3}\cdot \Phi_{p^2}^{k_4}\cdot \Phi_{p^2q}^{k_5}.$$
Let $p^*$ be the unique integer with $1 \leq p < q$
such that $p p^* \equiv 1 \pmod q$ and
$\mathcal{I}(f_k)$ be the integer interval:
\begin{itemize}
\item $[1,\ 1]$ for $k \in \{0, 2, 4, 14, 18, 62 \,\}$,
\item $[0,\ 1]$ for $k \in \{10, 12, 16, 26, 44, 46, 48, 50, 52, 58, 60 \,\}$,
\item $[-2,\ 2]$ for $k \in \{ 9, 35, 39, 41, 43, 57 \,\}$,
\item $[1,\ \min(p,\ q)]$ for $k \in \{ 6, 30 \,\}$,
\item $[0,\ \min(p,\ q)]$ for $k \in \{ 28, 54 \,\}$,
\item $[\min(\lfloor q / p \rfloor,\ 1), \min(\lfloor (q - 1) / p \rfloor + 1,\ p)]$ for $k = 20$,
\item $[1,\ \min(p^2,\ q)]$ for $k = 22$,
\item $[-\min(p,\ q - p^*),\ \min(p,\ p^*))]$ for $k = 24$,
\item $[-\min(p,\ p^*) - \min(p,\ q - p^*),\ \min(p,\ p^*) + \min(p,\ q - p^*)]$ for $k = 25$,
\item $[-\min(p,\ q),\ \min(p, q)]$ for $k = 29$,
\item $[-\min(p,\ q,\ q ({\rm mod~}p^2),\ p^2-q ({\rm mod~}p^2)),\ \min(p,\ q)]$ for $k = 38$,
\item $[-1,\ 1]$ otherwise.
\end{itemize}
Then $\mathcal{C}_0(f_k) = \mathcal{I}(f_k)$ except for $k = 38$ and $q = 2$.
If $q = 2$, $\mathcal{C}_0(f_{38}) = \mathcal{C}(f_{38}) = \{ -2, 0, 1, 2 \,\}$.
\medskip
We have $\mathcal{C}(f_k) = \mathcal{C}_0(f_k)$ except for the following
cases (where $\mathcal{C}(f_k) = \mathcal{C}_0(f_k) \setminus \{0\}$):
\begin{itemize}
\item $k = 1$,
\item $k \in \{13, 17, 29, 33, 61 \,\}$ and $p = 2$,
\item $k \in \{11, 40, 43, 59 \,\}$ and $q = 2$,
\item $k \in \{8, 24, 34 \,\}$ and $\min(p,\ q) = 2$,
\item $k = 9$ and $\min(p,\ q) \leq 3$,
\item $k = 25$ and $p \leq 3$ and $q \neq 2$,
\item $k = 38$ and $p = 2$ and $q = 3$,
\item $k = 41$ and $q \leq 3$.
\end{itemize}
\end{Thm}

\subsection{Earlier work on $x^{p^2q}-1$}
The only earlier work we are aware of is that by Kaplan \cite{K2},
who proved that if $p\neq q$, then $B(p^2q)=\min(p^2,q)$. He first remarks
that since $B(pq)=\min(p,q)$, it remains to deal with the 48 divisors of
$x^{p^2q}-1$ that do not divide $x^{pq}-1$. For those in his Table 1 he gives
an upper bound on the height. Since as we have seen in the proof of Theorem \ref{MAIN},
$H(f_{j+48})=H(f_j)$ for $0\le j\le 15$, it is actually enough to deal with only 32 divisors
(namely those in our Table 1B and 1C). A further remark is that where in his Table 1, $p$ is
given as upper bound, one needs $\min(p,q)$ (as not always $p\le \min(p^2,q)$). As Kaplan
pointed out to the authors it is easy to see that this replacement can be made. On doing so
and comparing with our results the upper bound he gives for the height are seen to be
equalities, except (for certain choices of $p$ and $q$) in the cases listed in Table 2.\\
%\vfil\eject

\centerline{{\bf Table 2}}
\begin{center}
\begin{tabular}{|c|c|c|}
\hline
$f$ & $H(f)$ & $H(f)$\\
\hline
& {\rm Kaplan} & {\rm exact}\\
\hline
$20$ & $\le \min(p,q)$ & $\min(p,[{q-1\over p}]+1)$\\
\hline
$24$ & $\le \min(p,q)$ & $\min(p,p^*)$\\
\hline
$25$ & $\le \min(p,q)$ & $\min(p,p^*)+\min(p,q-p^*)$\\
\hline
$45$ & $\le 2$ & 1\\
\hline
\end{tabular}
\end{center}
Note that two of the three `challenging' polynomials mentioned in the introduction do not
appear in the table. For $f_{38}$ it is easy to see that $H(f_{38})=\min(p,q)$ (but
challenging to determine ${\cal C}(f_{38})$). For $f_{43}$ it is easy to see that
$H(f_{43})\le 2$, but challenging to establish that $H(f_{43})=2$.
Of course in order to compute $B(p^2q)$ it is not the best strategy to compute
$H(f)$ exactly for every divisor of $x^{p^2q}-1$.

\section{Heights of divisors of $x^n-1$}
\label{height2010}
For a polynomial $f\in \mathbb Z[x]$, we define
$$H^*(f)=\max\{H(g):g|f {\rm ~and~}g\in \mathbb Z[x]\}.$$
Put $B(n)=H^*(x^n-1)$. So far little
is known about this function. Pomerance and Ryan \cite{PR} have established
the following three results concerning $B(n)$, the fourth is due to Justin \cite{Justin} and,
independently, Felsch and Schmidt \cite{FS}.
\begin{Thm}$~$\\
{\rm 1)} {\rm (\cite{PR}).} Let $p<q$ be primes. Then $B(pq)=p$.\\
{\rm 2)} {\rm (\cite{PR}).} We have $B(n)=1$ if and only if $n=p^k$.\\
{\rm 3)} {\rm (\cite{PR}).} We have
$$\lim_{n\rightarrow \infty}\sup {\log \log B(n)\over \log n/\log \log n}=\log
3.$$
{\rm 4)} {\rm (\cite{FS, Justin}).} $B(n)$ is bounded above by a function that does not depend on the largest prime factor of
$n$.\\
{\rm 5)} {\rm (\cite{K2}).} Let $p$ and $q$ be different primes. Then $B(p^2q)=\min(p^2,q)$.
\end{Thm} 
In their paper Pomerance and Ryan observe that from their limited numerical data it
seems that part 5 holds. This was subsequently proven by Kaplan \cite{K2}. Our work presented
here leads to a reproof. Kaplan's paper contains various further results on $B(n)$.\\
\indent For a polynomial $f\in \mathbb Z[x]$, we define
$$H^*_{\pm}(f)=\max\{|H_{\pm}(g)|:g|f {\rm ~and~}g\in \mathbb Z[x]\}.$$
Furthermore we define $B_{\pm}(n)=H^*_{\pm}(x^n-1)$. Numerical observations
suggest that often $B_{+}(n)>B_{-}(n)$, and this is our main motivation for
introducing these functions. In fact, if $p<q$ are primes, then
$B_{+}(pq)=p$ and $B_{-}(pq)=2$.

\section{Flat divisors of $x^n-1$}
The present article suggests that many divisors of $x^n-1$ are flat. It seems therefore
natural to try to obtain an estimate for the
number of flat divisors of $x^n-1$.\\
\indent The following result offers a modest contribution in this direction.
\begin{Thm}
\label{blaat}
Let $p$ and $q$ be distinct primes.
Let $f_e$ be the number of flat monic divisors of $x^{p^{e}q}-1$. Then
$f_{e+1}\ge 2f_{e}+2^{e+2}-1$.
\end{Thm}
{\it Proof}. Let $f(x)$ be a monic divisor of $x^{p^{e}q}-1$. Every divisor of $x^{p^{e+1}q}-1$ is of the form\\
a) $f(x)$;\\
b) $f(x)\Phi_{p^{e+1}}(x)$;\\
c) $f(x)\Phi_{p^{e+1}q}(x)$;\\ 
d) $f(x)\Phi_{p^{e+1}}(x)\Phi_{p^{e+1}q}(x)$.\\

\noindent Lower bounds for the number of flat divisors
 amongst the types a,b,c and d are considered below:\\
a) The divisors of this form contribute $f_e$ to $f_{e+1}$.\\
b) Note that we can write $f(x)\Phi_{p^{e+1}}(x)=f(x)\Phi_p(x^{p^e})$. Suppose $f(x)$
divides $x^{p^e}-1$. Since $f(x)\Phi_p(x^{p^e})|x^{p^{e+1}}-1$ it is flat by Theorem \ref{riant}.
Since $x^{p^e}-1$ has $2^{1+e}$ monic divisors, we
see that there are at least $2^{1+e}$ flat divisors of  $x^{p^{e+1}q}-1$ of the
form b.\\
c) Note that we can write $f(x)\Phi_{p^{e+1}q}(x)=f(x)\Phi_{pq}(x^{p^e})$. Suppose $f(x)$
divides $x^{p^e}-1$. In case $f(x)=x^{p^e}-1$, then $H(f(x)\Phi_{pq}(x^{p^e}))=2$ by 
Lemma \ref{nulgat}. In
the remaining case deg$(f)<p^e$ and by (\ref{bier}) and Theorem \ref{riant} we infer that
$H(f(x)\Phi_{pq}(x^{p^e}))=H(f(x))=1$. We conclude that there are at least $2^{1+e}-1$ flat divisors of  $x^{p^{e+1}q}-1$ of the
form c.\\
d) We have $\Phi_{p^{e+1}}(x)\Phi_{p^{e+1}q}(x)=(x^{p^{e+1}q}-1)/(x^{p^eq}-1)$.
In case $f(x)=x^{p^eq}-1$, then $H(f(x)\Phi_{p^{e+1}}(x)\Phi_{p^{e+1}q}(x))=H(x^{p^{e+1}q}-1)=1$.
In the remaining cases we find, by (\ref{bier}), that 
$H(f(x)\Phi_{p^{e+1}}(x)\Phi_{p^{e+1}q}(x))=H(f(x))$. Thus there are at least $f_e$
flat divisors of $x^{p^{e+1}q}-1$ of the form d.\\
\indent On adding the contributions of each of the four types a,b,c and d, the result follows. \qed\\

\noindent {\tt Remark 1}. The above argument with $e=1$ 
in combination with Theorem \ref{plakkerig} leads to the following list of 35
divisors of $x^{p^2q}-1$ that are flat:\\
a) $f_0,\ldots,f_{15}$, excluding $f_6$ and $f_9$\\
b) $f_{16},f_{17},f_{18},f_{19}$\\
c) $f_{32},f_{33},f_{34}$\\
d) $f_{48},\ldots,f_{63}$, excluding $f_{54}$ and $f_{57}$.\\
Note that the full list of flat divisors is longer.\\

\noindent {\tt Remark 2}. By induction one easily proves that for $e\ge 2$
we have $$f_e\ge 2^{e-1}f_1+(4e-5)2^{e-1}+1.$$ By Theorem \ref{plakkerig} we have 
$f_1=14$ and hence it follows that $f_e\ge (4e+9)2^{e-1}+1$. The total
number of divisors of 
$x^{p^eq}-1$ is $2^{2+2e}$, denote this by $n_e$. Then $f_e\gg \sqrt{n_e}\log n_e$. 
Can one improve on this ?

\section{A variation}
We have $H(f_6)=\min(p,q)=B(pq)$. Likewise we have $H(f_{22})=\min(p^2,q)=B(p^2q)$. Both $f_6$ and $f_{22}$
are special in the sense that they have only non-negative coefficients. It might therefore be more reasonable to
consider only balanced divisors of $x^n-1$, that is divisors having both positive and negative coefficients. 
Let us denote this
analogue of $B(n)$ by $B'(n)$. Put
$$C(n)=\max\{|{\cal C}_0(f)|-1: f|x^n-1,~f{\rm ~is~balanced}\}.$$
\begin{Thm}
We have\\ 
{\rm 1)} $B'(pq)=2$ and $C(pq)=4$.\\ 
{\rm 2)} $B'(p^2q)=B_{-}(p^2q)=\min(p,p^*)+\min(p,q-p^*)$ and $C(p^2q)=2B'(p^2q)$.
\end{Thm}
This result is a consequence of the inequality $\min(p,p^*)+\min(p,q-p^*)\ge \min(p,q)$ and Theorem \ref{MAIN}.
It does not follow from earlier work in this area (\cite{K2, PR,ryanryan}).\\

\noindent {\tt Acknowledgement}. The bulk of this paper was written in August/September 2008 during an internship of the
first author with the second author. The initial aim was to prove the conjecture of Ryan and Pomerance
that $B(p^2q)=\min(p^2,q)$. This was relatively soon achieved, but then we learned that independently this
had already been done by Kaplan \cite{K2}. Then the aim became to compute the maximum and minimum coefficient
of each of the 64 monic divisors of $x^{p^2q}-1$. With the more recent focus on coefficient convexity 
(see, e.g., \cite{B, BZ, GM, Rosset}) in mind the aim
was even set higher: to compute the coefficient set of all of the divisors.\\
\indent We like to thank the interns Richard Cartwright (2008) and Oana-Maria Camburu (2010) for helpful
remarks. However, our greatest indebtedness is to Yves Gallot for his computational 
assistance. In particular, he numerically verified Theorem \ref{MAIN} in case $\max(p,q)<200$.
Merci beaucoup, Yves!

\medskip\noindent {\footnotesize Achtern Diek 32, D-49377 Vechta, Germany.\\
e-mail: {\tt andreasd@uni-bonn.de}\\

\medskip\noindent {\footnotesize Max-Planck-Institut f\"ur Mathematik,\\
Vivatsgasse 7, D-53111 Bonn, Germany.\\
e-mail: {\tt moree@mpim-bonn.mpg.de}\\

\end{document}